\documentclass[11pt]{article}

\usepackage[english]{babel}
\usepackage[T1]{fontenc}
\usepackage[utf8]{inputenc}
\usepackage{amsmath,amsfonts,amsthm,graphicx,color,mdwlist}
\usepackage{qsymbols,booktabs,psfrag,color}
\usepackage[%
            linktocpage,%
            bookmarks=true,%
            colorlinks=true,%
            ]{hyperref}
\usepackage{textcomp,xspace}

\def \eref#1{(\ref{#1})}
\def \1{{\bf 1}}
\def \l{\left}
\def \r{\right}

\def \ben{\begin{eqnarray}}
\def \een{\end{eqnarray}}
\def \be{\begin{eqnarray*}}
\def \ee{\end{eqnarray*}}

\def \beq{\begin{equation}}
\def \eq{\end{equation}}
\def \Per{{\sf Peri}}
\def \Ext{{\sf Ext}}
\def \Mult{{\sf Multinomial}}
\def \Pois{{\sf Poisson}}

\def \uniform{{\sf uniform}}
\def \Coeffs{\textsf{Coeffs}}
\def \CT{{\sous{\star}{\T}}}
\def \Circle{{\sf{Circle}}}
\def \BB{\mathcal{B}}
\def \CC{\mathcal{C}}

\def \Barany{B\'ar\'any }
\def \Voronoi{Vorono\"\i\ }

\def \app#1#2#3#4#5{\begin{array}{rccl} 
    #1:&#2&\longrightarrow&#3\\ &#4&\longmapsto&#5
  \end{array}}
\def \jfmatop#1#2{\genfrac{}{}{0pt}{}{#1}{#2}}
\def \Proof{{\noindent \bf Proof. \rm}}
\def \Proofof#1{\noindent \textbf{Proof of #1.}}
\newtheorem{lem}{Lemma}[section]
\newtheorem{prop}[lem]{Proposition}%[section]
\newtheorem{theo}[lem]{Theorem}%[section]
\newtheorem{cor}[lem]{Corollary}%[section]
\newtheorem{defi}[lem]{Definition}%[section]
\theoremstyle{definition}
\newtheorem{exm}[lem]{Example}
\newtheorem{rem}[lem]{Remark}%[section]
\def \MT{{\cal M}_{\cal T}}
\def \T{{\cal T}}

\def \dd{\xrightarrow[n]{(d)}}
\def \weak{\xrightarrow[n]{(weak)}}
\def \proba{\xrightarrow[n]{(proba.)}}
\def \as{\xrightarrow[n]{(a.s.)}} 
\def \iid{i.i.d.\@\xspace}
\def \rv{r.v.\@\xspace}
\def \AS{a.s.\@\xspace}
\newcommand{\Rel}{\mathbb{Z}}

\newcommand{\Rea}{\mathbb{R}}
\newcommand{\Com}{\mathbb{C}}

\def \Conv {\textsf{Conv}\xspace}
\def \BConv {\textsf{BConv}\xspace}

\def \sous#1#2{\mathrel{\mathop{\kern 0pt#1}\limits_{#2}}}
\def \sur#1#2{\mathrel{\mathop{\kern 0pt#1}\limits^{#2}}}
\def \DA{{\sf DA}}

\DeclareMathOperator{\Seg}{Seg}
\DeclareMathOperator{\Nei}{Nei}

\DeclareMathOperator{\cov}{cov}

\def \bB{{\bf B}}
\def \sK{{\sf K}}

\newqsymbol{`E}{\mathbb{E}}

 \newqsymbol{`P}{\mathbb{P}}

\newqsymbol{`R}{\mathbb{R}}

\newqsymbol{`e}{\varepsilon}
\newqsymbol{`o}{\omega}

\begin{document}
\begin{center}

\LARGE{\bf Compact convex sets of the plane\\[1mm]
  and probability theory}\medskip\normalsize
\[\begin{array}{ll}
\textrm{\Large Jean-Fran\c{c}ois Marckert}& 
\textrm{\Large ~~~~~~~~~~David Renault}\end{array}\]
\textrm{CNRS, LaBRI, Universit\'e de Bordeaux}\\
\textrm{351 cours de la Lib\'eration}\\
\textrm{33405 Talence cedex, France}\\
\textrm{email: name@labri.fr}
 \end{center}

\begin{abstract} 
  The Gauss-Minkowski correspondence in $`R^2$ states the existence of
  a homeomorphism between the probability measures $\mu$ on $[0,2\pi]$
  such that $\int_0^{2\pi} e^{ix}d\mu(x)=0$ and the compact convex
  sets (CCS) of the plane with perimeter~1. In this article, we bring
  out explicit formulas relating the border of a CCS to its
  probability measure.
  As a consequence, we show that some natural operations on CCS -- for
  example, the Minkowski sum -- have natural translations in terms of
  probability measure operations, and reciprocally, the convolution
  of measures translates into a new notion of convolution of CCS.
  Additionally, we give a proof that a polygonal curve associated with
  a sample of $n$ random variables (satisfying $\int_0^{2\pi}
  e^{ix}d\mu(x)=0$) converges to a CCS associated with $\mu$ at speed
  $\sqrt{n}$, a result much similar to the convergence of the
  empirical process in statistics. Finally, we employ this
  correspondence to present models of smooth random CCS and
  simulations. 
 \end{abstract}

\section{Introduction}

Convex sets are central in mathematics: they appear everywhere\,!
Nice overviews of the topic have been provided by
Busemann~\cite{bus}, P\'olya \cite{pol} and Pogorelov \cite{Pog}.
In probability theory, \textit{compact convex sets} (CCS) appear in
1865 with Sylvester's question \cite{Syl}: for $n=4$ points chosen
independently and at random in the unit square $K$, what is the
probability that these $n$ points are in convex position~? The
question can be generalised to various shapes~$K$, different values of
$n$, and other dimensions.
It has been recently solved by Valtr \cite{Valtr1,Valtr2}
when $K$ is a triangle or a parallelogram and by Marckert
\cite{MarckertCircle} when $K$ is a circle (see also
\Barany\cite{BAR}, Buchta \cite{Buchta} and
\Barany\cite{BaranyRandomPolytopes}).
Random CCS also show up as the cells of the \Voronoi diagram of a
Poisson point process (see Calka \cite{Calka}), and in the problem of
determining the distribution of convex polygonal lines subject to some
constraints. For example, when the vertices are constrained to belong
to a lattice, the problem has been widely investigated (Sinai
\cite{Sinai}, \Barany \& Vershik \cite{B-V}, Vershik \&
Zeitouni~\cite{V-Z}, Bogachev \& Zarbaliev~\cite{Bo-Za}). 
Another combinatorial model related to this question is based on the
\textit{digitally convex polyominos} (DCPs). The DCP associated to a
convex planar set~$C$ is the maximal convex polyomino with vertices in
$\mathbb{Z}^2$ included in $C$. Let $D_n$ be the set of DCPs with
perimeter $2n$.  In a recent paper, Bodini, Duchon \& Jacquot
\cite{BDJ} investigate the limit shape of uniform DCPs taken in $D_n$
under the uniform distribution~$\mathbb{U}_n$. Even if not convex,
these polyominos can be seen as discretisation of~CCS.

\medskip

All these models possess the same drawbacks: they are discrete models
(polygonal, except for DCP) and their limit when the size parameter
goes to $+\infty$ are deterministic shapes.
To our knowledge, no model of random non-polygonal CCS have been
investigated yet. One of the goals of this article is to develop tools
that allow one to provide examples of such models, and this goal is
attained in the following manner~:

\medskip

\noindent\textbullet~~First, we state a connection between the CCS of
the plane and probability measures. Theorem~\ref{th:mes-conv} asserts
that the set of CCS of the plane having perimeter 1, considered up to
translation, is in one-to-one correspondence with the set $\MT^0$ of
probability distributions $\mu$ on the circle $`R/(2\pi \mathbb{Z})$
satisfying $\int_{0}^{2\pi}\exp(ix)d\mu(x)=0$. This famous theorem,
revisited in Section \ref{sec:mes-ci}, is sometimes called in the
literature the {Gauss-Minkowski} Theorem (cf. Vershik~\cite{V-Z} and
Busemann \cite[Section~8]{bus}), and the measure $\mu$ is called the
surface area measure of the CCS \cite{MOS}.
Moreover, the bijection is an homeomorphism when both sets are
equipped with natural topologies.  In this article, we provide an
explicit parametrisation of a CCS in terms of the distribution
function of $\mu$. This perspective brings out a new and important
relation between the CCS with perimeter $1$ and probability measures,
differing in this from the more generic ``arbitrary total mass''
measures.

\medskip

\noindent\textbullet~~This connection with probability theory 
appears therefore as a natural tool to define new operations on CCS
and revisit numerous known results that were proved using geometrical
arguments.
For instance, the set $\MT^0$ is stable by convolution and mixture.
This induces natural operations on CCS that one may also qualify of
\textit{convolution} and \textit{mixture}. As a matter of fact, the
mixture of CCS defined in this way coincides with the Minkowski
addition (Section \ref{sec:Mixt}), and Minkowski symmetrisation
simply maps a CCS associated to a measure $\mu$ onto the CCS
associated with $\frac{1}{2}(\mu+\mu(2\pi-\,.\,)$ (Proposition
\ref{pro:mix=Min}).
The notions of \textit{convolution of CCS} and \textit{symmetrisation
  by convolution} (Sections \ref{sec:Conv} and \ref{sec:symm}) appear
to be new and provide a new proof of the isoperimetric inequality
(Theorem \ref{thm:Mixture2}). Roughly, the CCS obtained by convolution
of two CCS has a radius of curvature function equal to the convolution
of the curvature functions of these two CCS.

\medskip

\noindent\textbullet~~The probabilistic approach also
allows one to prove stochastic convergence theorems for models that
differ radically from the ones mentioned earlier.
Consider for instance $\mu \in \MT^0$, and take~$n$ random variables
$\{X_j, j=0,\dots,n-1\}$ \iid according to~$\mu$. Let
$\{\widehat{X}_j,j=0,\dots,n-1\}$ be the $X_k$'s reordered in
$[0,2\pi)$. Let $B_n$ be the curve formed by the concatenation of the
vectors $e^{i\widehat{X}_j}$. We show that the curve $B_n$ rescaled by
$n$ converges when $n \to\infty$ to the boundary $\BB_\mu$ of a CCS
associated with $\mu$ (Theorem \ref{thm:cv-emp} and Corollary
\ref{cor:cv-emp}). This convergence holds at speed $\sqrt{n}$ and has
Gaussian fluctuations (Theorem \ref{thm:cv-emp}).
As a generalisation, every distribution on $\mathbb{C}$ with mean~0
can be sent on a CCS by a second correspondence (which is not
bijective) (Section \ref{sec:mes-C}). Again, the appropriate point of
view consists in considering the boundary of the CCS as the limit of
the curve associated with a sample of $n$ random variables (\rv)
sorted according to their argument.

\medskip

\noindent\textbullet~~The 
last part of this paper (Section \ref{sec:mrc}) is devoted to the
investigation of models of random CCS that stem from the aforesaid
connection. Our first model is a model of random polygons defined as
follows: take $\{z_j, j=0,\dots,n-1\}$ \iid according to a
distribution $\nu$ in $\mathbb{C}$. Let $\{y_i=z_{i+1 \bmod{n}}-z_i,
i=0,\dots,n-1\}$ and $\{\widehat{y}_j,j=0,\dots,n-1\}$ the $y_i$'s
sorted according to their argument. The $\widehat{y}_i's$ are the
consecutive vector sides of the polygonal CCS with vertices
$\{\sum_{j=0}^{d}\widehat{y}_j,d = 0,\dots,n-1\}$.  When $n\to
\infty$, a rescaled version of this CCS converges in distribution to a
deterministic CCS (Theorems \ref{thm:cv-conve} and
\ref{thm:conv-conv}). We discuss the finite case in
Section~\ref{sec:cccc}.

\noindent\textbullet~~Another model
results from the role that Fourier series play in the representation
of the boundaries of CCS. For a \rv $X$ with values in $[0,2\pi]$ and
distribution $\mu$, the Fourier coefficients of $\mu$, namely
$a_n(\mu)=`E(\cos(nX))$ and $b_n(\mu)=`E(\sin(nX))$, are well defined
for any $n\geq 0$. Our bijection between CCS and measures in hand, the
question of designing a model of random CCS is equivalent to that of
designing a model of random measure $\mu$ satisfying \AS
$\int_{0}^{2\pi}\exp(ix)d\mu(x)=0$ (equivalently $a_1(\mu) = b_1(\mu)
= 0$ \AS). Nevertheless to design a model of random measures $\mu$
satisfying these constraints is not equivalent to design random
Fourier coefficients $(a_n,b_n,n\geq 0)$ since these latter may not
correspond to those of a probability measure. In
Section~\ref{sec:mrc}, we explain how this can be handled, and provide
several models of random CCS that are not random polygons. \medskip

\noindent \textbf{Notations.}~
``CCS'' will always be used for ``compact convex set of the plane
$`R^2$''.  
We assume that all the mentioned \rv are defined on a
common probability space $(\Omega,{\cal A},`P)$, and denote by $`E$
the expectation. For any probability distribution~$\mu$, $X_\mu$
designates a \rv with distribution $\mu$. We write $X\sim \mu$ to say
that $X$ has distribution $\mu$. The notations $\dd, \proba,\weak$
stand for the convergence in distribution, in probability, and the
weak convergence.

\section{Correspondence between CCS and distributions}

We start this section by recalling some simple facts concerning CCS
and measures on the circle $`R/(2\pi \mathbb{Z})$.  Thereafter we
state the Gauss-Minkowski theorem (Theorem \ref{th:mes-conv}) which
establishes a correspondence between measures and CCS, and we provide
a new proof based on probabilistic arguments. In Section
\ref{sec:fourier} we express the area of a CCS thanks to the Fourier
coefficients of the associated measure. Finally in Section
\ref{sec:conver-conv} we state one of the main results of the paper
(Theorem~\ref{thm:cv-emp}): under some mild hypotheses, it ensures the
convergence of the trajectory made of $n$ \iid increments sorted
according to their arguments and rescaled by $n$ to a limit CCS
boundary at speed~$\sqrt{n}$.

\subsection{ CCS of the plane}
A subset $S$ of $`R^2$ is a \textit{convex} set if for any $z_1,z_2\in
S$, the segment $[z_1,z_2]\subset S$.  In this paper, we are
interested only in CCS of the Euclidean plane $`R^2$.  Let $\Seg$ be
the set of bounded closed segments, and $\Nei$ be the set of
CCS with \underline{n}on \underline{e}mpty \underline{i}nteriors. The
union $\Seg \cup \Nei$ forms the set of all CCS of $`R^2$. \par

% 
% Nei case
%
For $S \in \Nei$, $S^\circ$ will designate the interior of $S$,
and $\partial S={S}\setminus S^\circ$ the boundary of~$S$. We call
\textit{parametrisation} of $\partial S$, a map $\gamma:[a,b]\to
\partial S$ for some interval $[a,b]\subset \mathbb{R}$, such that
$\gamma(a)=\gamma(b)$ and such that $\gamma$ is injective from $[a,b)$
to $\partial S$.
The length of $\partial S$ is well defined, finite and positive, and
is called the \textit{perimeter} of $S$ and denoted~$\Per(S)$.  
It may be used to provide a \textit{natural parametrisation} of
$\partial S$, that is to say a function $\gamma:[0,|\partial S|] \to
\partial S$, continuous and injective on $[0,|\partial S|]$, such that
$\gamma(0)=\gamma(|\partial S|)$ and such that the length of
$\{\gamma(t), t\in[0,s]\}$ is equal to $s$ for any $s\in[0, |\partial
S|]$.
%
% Seg case
%
For $S \in \Seg$, the notion of natural parametrisation also
exists, but it is different. For technical reasons, we choose the
following one: The natural parametrisation of a segment $[a,b]$ is
defined to be $\gamma(t)=a(1-\frac{t}{|b-a|})+b\frac{t}{|b-a|}$ on
$[0,|b-a|]$ and $\gamma(t)=a(\frac{t}{|b-a|}-1)+b(2-\frac{t}{|b-a|})$
on $[|b-a|,2|b-a|]$, as if the segments were thick and two-sided. 
In this case, we define $\Per(S) = 2|b-a|$.

\begin{defi} The \textit{boundary} $B$ of $C \in \Nei$ is defined as
  $B=C\setminus{C^\circ}$. The boundary $B$ of $C=[a,b] \in \Seg$ is
  $C$ itself. 
\end{defi}

The boundary of a CCS is equal to the path induced by its natural
parametrisation, and its perimeter is the length of this path.

\subsection{Measures on the circle}
\label{sec:mes-ci}
Let ${\T}$ be the circle $`R/(2\pi\mathbb{Z})$ equipped with the
quotient topology, and $\MT$ be the set of probability distributions
on $\T$. 
The weak convergence on $\MT$ is defined as usual: $(\mu_n,n\geq 0)\weak \mu$ in $\MT$ if for any 
bounded continuous function $f:\T\to `R$, $\int_\T f d\mu_n\to  \int_\T f d\mu$.
Let  $\mu\in \MT$, and consider 
\[\app{F_\mu}{\T}{[0,1]}{x}{\mu([0,x])}\] 
be the \textit{cumulative distribution function} (CDF) of $\mu$. Let $\mathcal{I}_\mu$ be the set of points 
of continuity of $F_\mu$, where by convention, $0\in \mathcal{I}_\mu$ if $F_{\mu}(0)=\mu(\{0\})=0$. If $\mu_n\weak\mu$ in
$\MT$, then it can not be deduced that $F_{\mu_n}\to F_\mu$ pointwise on $\mathcal{I}_\mu$ since  $\delta_{2\pi}=\delta_0$ in $\MT$. 
What is still true, is that
\[F_{\mu_n}(y)-F_{\mu_n}(x)\to F_{\mu}(y)-F_{\mu}(x), \textrm{ for
  any }(x,y)\in \mathcal{I}_\mu.\] 
A function $F:[0,2\pi)\to `R$ is a CDF of some distribution $\mu\in \MT$ if it is right continuous, non decreasing on $[0,2\pi]$, satisfies
$0\leq F(0)\leq 1$, $F(2\pi-)=1$ 
(see Wilms \cite[p.4-5]{WilmsFractional} for additional information and references).

Consider the function
\begin{equation}\label{eq:conv-measure2} 
\app{Z_\mu}{[0,1]}{\mathbb{C}}{t}{Z_\mu(t) = \displaystyle\int_0^t \exp(iF_\mu^{-1}(u)) du,}
\end{equation}
where $F_\mu^{-1}$ is the \textit{standard generalised inverse} of $F_\mu$:
\[\app{F_\mu^{-1}}{[0,1]}{[0,2\pi)}{y}{F_\mu^{-1}(y) :=\inf \{x \geq
  0~:~F_\mu(x)\geq y\}.}\]
The range $\BB_\mu$ of $Z_\mu$ is the central object here:
\[\BB_\mu:=\l\{ Z_\mu(t),~t\in[0,1]\r\}.\]
Since the modulus of $Z_\mu'$ is 1, $Z_\mu$ is the natural parametrisation of $\BB_\mu$ 
and $\BB_\mu$ has length~1. \par

Let $\Conv $ be the set of CCS of the plane containing the origin,
lying above the $x$-axis, and whose intersection with the $x$-axis is
included in $`R^+$.  Denote by $\Conv(1)$ the subset of $\Conv$ of CCS
having perimeter 1, and by $\BConv$ the set of their corresponding
boundaries. Set
\[\MT^0=\l\{\mu \in\mathcal{M}[0,2\pi]~, \int_{0}^{2\pi-} \exp(i\theta)dF_\mu(\theta)=0\r\}\]
the subset of $\MT$ of measures having Fourier transform equal to 0 at
time 1.

\subsection{Probability measures and CCS}

Probability distributions on $`R$ are characterised by their Fourier
transform, and convergence of Fourier transforms characterises weak
convergence by the famous Lévy's continuity Theorem. The following
Theorem gives a similar characterisation of measures in $\MT^0$ by
their representation as CCS of the plane.

\begin{theo} \label{th:mes-conv}
1)  The map
\[\app{\BB}{\MT^0}{\BConv(1)}{\mu}{\BB_\mu}\] 
is a bijection. \\
2) $\BB$ is an homeomorphism from $\MT^0$ (equipped with
  the weak convergence topology) to $\BConv(1)$ (equipped with the
  Hausdorff topology on compact sets). \\
3)  The function $\Gamma$ from $\Conv(1)$ to $\BConv(1)$ which sends a CCS to its boundary is an homeomorphism for the Hausdorff topology, and then 
\[\app{\CC}{\MT^0}{\Conv(1)}{\mu}{\CC_\mu:=\Gamma^{-1}(\BB_\mu)}\]
is an homeomorphism.
 \end{theo}

This theorem sometimes called ``Gauss-Minkowski'' in the literature
can be found in a slightly different form in Busemann \cite[Section
  8]{bus}. 
The integral formula \eref{eq:conv-measure2} giving the
parametrisation of the CCS in terms of $F_\mu^{-1}$, which is
central here, seems to be new. We provide a proof of Theorem
\ref{th:mes-conv} in probabilistic terms at the end of this section.

In Busemann, this theorem is stated more generally in $\Rea^n$,
where the measures range over the unit sphere of $\Rea^n$ and verify a
set of properties, which in $\Rea^2$ sum up to $\int_0^{2\pi}
e^{ix}d\mu(x) = 0$. 
The measure $\mu$ is called the surface area measure \cite{MOS} of the CCS $\CC_\mu$, and is defined for more general convex sets in any dimension. 

\begin{rem}
The map $\BB$ that one may see as a ``curve'' transform, may be
extended to $\mathcal{M}[0,2\pi]$, the set of measures on $[0,2\pi]$;
in this case $\BB(\mathcal{M}[0,2\pi])$ is the set of continuous
almost everywhere differentiable curves of length $1$, starting at
the origin, having a positive argument in a neighbourhood of 0, and where
along an injective parametrisation, the argument of the tangent is
non decreasing\footnote{The Fourier transform $t\mapsto \Psi_\mu(t)$ also
  defines a curve $\{\Psi_\mu(t): t \in A\}$ in the plane, for any
  interval $A$. This curve is different from $\CC_\mu$, for any $A$.}.
\end{rem}

There exists another formula for $Z_\mu$ in terms of expectations of
\rv, that we will use as a guideline throughout the paper. 
Recall that if $U\sim \uniform[0,1]$ then $F_{\mu}^{-1}(U)\sim \mu$, and then
\beq\label{eq:Change_variable}
 Z_\mu(t)=`E\l( \1_{U\leq t} \exp(iF_\mu^{-1}(U))\r).
\eq
Since $x\leq F_\mu(y)$ is equivalent to $F_\mu^{-1}(x)\leq y$, we obtain that
\be
Z_\mu(F_\mu(t))&=&`E\l( \1_{U\leq F_\mu(t)} \exp(iF_\mu^{-1}(U))\r)\\
               &=&`E\l( \1_{F^{-1}_\mu(U)\leq t} \exp(iF_\mu^{-1}(U))\r)\\  
               &=&`E\l( \1_{X_\mu\leq t} \exp(iX_\mu )\r). 
\ee 
The function $t\mapsto Z_\mu(F_\mu(t))$ plays an important role since
it encodes the extremal points of $\BB_\mu$ (see below). The function
$Z_\mu$ is somehow less pleasant since it can not be written directly
in term of $X_\mu$ on $[0,1]$. To see this, let 
\be
I_\mu=\{t \in [0;2\pi) ~\textrm{such that}~ \{u, u < t\}= \{F_\mu^{-1}(u)
< F_\mu^{-1}(t)\} \}
\ee
This corresponds to the set of $t$ where $F_\mu^{-1}(t)>F_\mu^{-1}(t-h)$
for any $h>0$ (or $t=0$). It can be shown that $I_\mu=\{F(t), t\in
[0,2\pi]\}$. Noticing that one can replace $\1_{U \leq t}$ by
$\1_{U<t}$ in \eref{eq:Change_variable}, we have
\ben\label{eq:for2}
Z_\mu(t)=`E\l( \1_{X_\mu < F^{-1}_\mu(t)} \exp(iX_\mu)\r) \textrm{ for }~ t \in I_\mu,
\een
Now we can characterise $\Ext(C)$ the set of extremal points of $C$.
\begin{lem}\label{lem:ext-points} For any $\mu\in \MT^0$, 
$\Ext(\CC_\mu)={\l\{Z_\mu(F_\mu(t)), t \in [0,2\pi]\r\}}$.
\end{lem}

\Proof 
From \eref {eq:Change_variable}, we see that $Z_\mu$ is linear on every interval inside the complement of $I_\mu$ in $[0,1]$: if $(t_1,t_2)$ is such an interval, for any $t\in[t_1,t_2]$,
\[Z_\mu(t)=Z_\mu(t_1)+(t_2-t)\frac{Z_\mu(t_2)-Z_\mu(t_1)}{t_2-t_1}.\]
Therefore, the points in the complement of $I_\mu$ are not extremal, 
and reciprocally, every non-extremal point lies on a segment inside $\BB_\mu$ and necessarily
belongs to the complement of $I_\mu$. Therefore
$\Ext(\CC_\mu)$ is equal to the closed set $\{Z_\mu(F_\mu(t), t \in [0,2\pi]\}$. ~$\Box$
\medskip

The curvature $k_\mu(t)$ of $\CC_\mu$ at time $t$, is given by
$\frac{1}{F'_\mu(F^{-1}_\mu(t))}$ when $F_\mu$ admits a derivative at
$F^{-1}_\mu(t)$; in particular, this means that when $\mu$ admits a
density $f_\mu$, then
$k_\mu(F_\mu(\theta))=1/f_\mu(F^{-1}_\mu(F_\mu(\theta)))=1/f_\mu(\theta)$,
which corresponds to the curvature at the point whose tangent has
direction $\theta$.

The real and imaginary parts $x_\mu(t)=\Re(Z_\mu(t))$ and
$y_\mu(t)=\Im(Z_\mu(t))$ of $Z_\mu(t)$ satisfy
\begin{equation}\label{eq:curve-x-y}
\left\{
\begin{array}{cl}
x_\mu(t) &= \int_0^t \cos(F^{-1}_\mu(u)) du  = \int_0^{F^{-1}_\mu(t)} \cos(v)
dF(v)   \\
y_\mu(t) &= \int_0^t \sin(F^{-1}_\mu(u)) du  = \int_0^{F^{-1}_\mu(t)} \sin(v)
dF_\mu(v)  .
\end{array}\right.
\end{equation}
the second equality in each line being valid only for $t\in I_\mu$.
\begin{figure}[ht]
\begin{center}
\includegraphics[height=5.5cm]{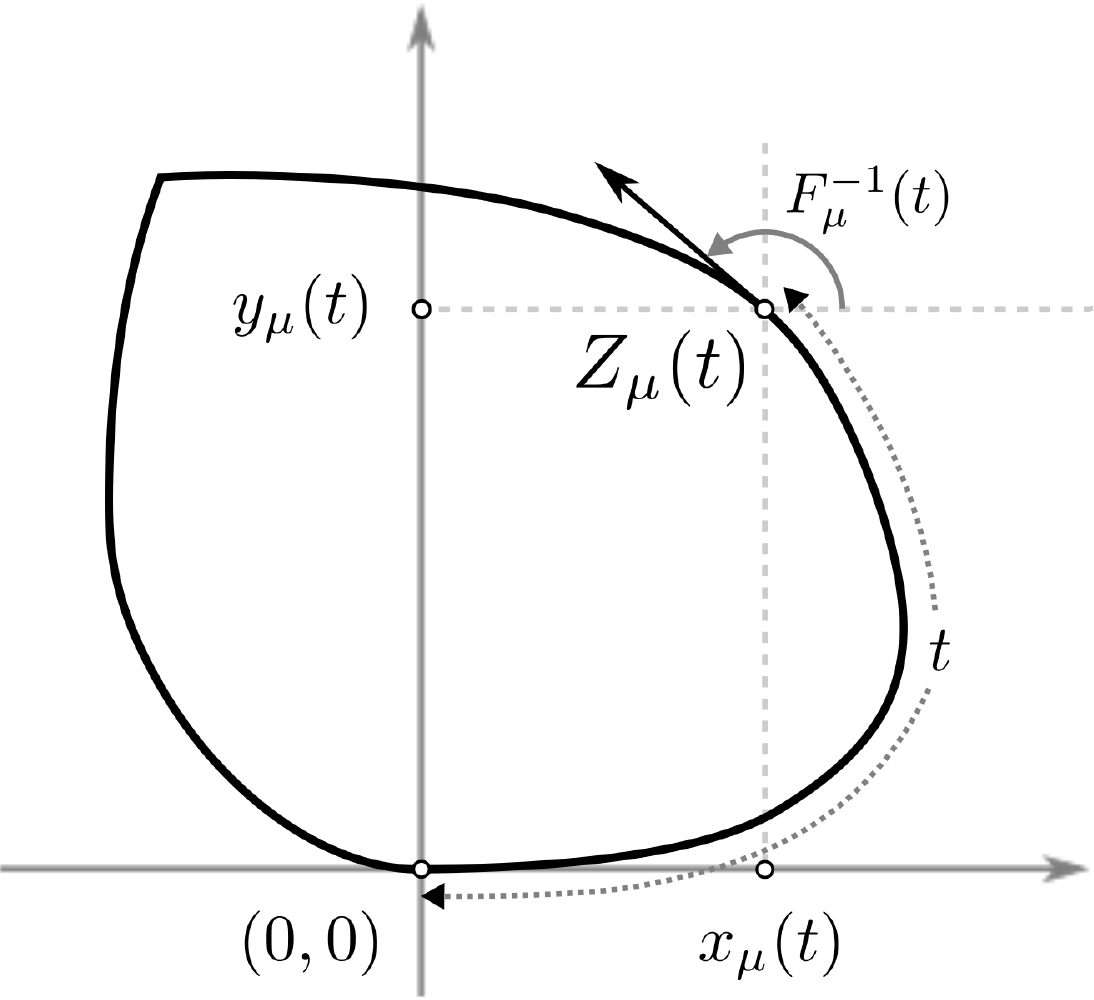}
\caption{A CCS $\CC_\mu$ for some measure $\mu$, $t$ gives the length
  of the curve $\BB_\mu$ between 0 and $Z_\mu(t)$ (in the
  trigonometric order), $F^{-1}_\mu(t)$ is then the direction of the
  tangent at time $t$.}
\label{fig:code-conv2}
\end{center}
\end{figure}

\Proofof{Theorem \ref{th:mes-conv} 1)}
The proof of $3)$ is immediate. We establish $1)$.

a) First, we prove that for any $\mu\in\MT^0$, $\BB_\mu$ is the
boundary of a CCS $\CC_\mu \in \Conv(1)$. A \textit{support
  half-plane} of $\BB_\mu$ is a half-plane $H$ intersecting $\BB_\mu$
on its border and such that $\BB_\mu \subset H$. The function $Z_\mu$
is continuous, and a simple analysis shows that $y_\mu$ is such that
$y_\mu(0)=y_\mu(1)=0$, and is increasing then decreasing over
$[0,1]$. Therefore, $\BB_\mu$ lies on the half plane above the
$x$-axis, which is a support half-plane of $\BB_\mu$. More generally,
for any $\theta\in [0,2\pi)$, $\mu_\theta(.)=\mu(.-\theta \mod 2\pi)$
  is still in $\MT^0$, and $\BB_{\mu_\theta}$ lies on the half plane
  above the $x$-axis. Therefore, for all $t\in[0,1)$, the line $D_t$
    passing through $Z_\mu(t)$ making an angle $F_\mu^{-1}(t)$ with
    the origin, is the border of a support half-plane of $\BB_\mu$.
Since $F^{-1}_\mu$ is right-continuous, $\BB_\mu$ is even tangent to
$D_t$.

We now show that $\BB_\mu$ is a simple curve or a segment: let $z$ be
such that $z=Z_\mu(t_1)=Z_\mu(t_2)$, for $t_1<t_2$. Then, by
definition \eref{eq:conv-measure2},
$\int_{[t_1,t_2]}\exp(iF_\mu^{-1}(u)) du =
\int_{[0,t_1]\cup[t_2,1]}\exp(iF_\mu^{-1}(u)) = 0$. Each of these
integrals is the weighted barycentre of a portion of the circle, both
portions being disjoint except at their extremities $t_1$ and
$t_2$. Since both barycentres are equal (to 0), the support of $\mu$
must be included in $\{t_1,t_2\}$. This implies that $F^{-1}_\mu(t_2)
= \pi + F^{-1}_\mu(t_1)$ and $\mu(\{t_2\})=\mu(\{t_1\})=1/2$. In other
words, the CCS is a segment of length $1/2$.  Therefore, when
$\BB_\mu$ is not a segment, it is a bounded Jordan curve that encloses
a bounded connected subset~$\CC_{\mu}$.
In this last case, $\BB_\mu$ is the border of $\CC_\mu$ and every
point of the border possesses a support half-plane, therefore
$\CC_\mu$ is convex (see for example 3.3.6 in \cite{MOS}).

b) The injectivity of $\BB$ is clear since if
$F_\mu^{-1}(t)=F_\nu^{-1}(t)$ for all $t\in[0,1]$, then
$\mu=\nu$. Now, let $B$ be a CCS boundary in $\BConv(1)$ and consider
the unique natural parametrisation $Z$ of $B$ in the counterclockwise
direction such that $Z(0)=Z(1)=0$. The map $Z$ has almost everywhere a
derivative $g$, and since it is continuous, $g$ is the derivative of
$Z$ in the distribution sense: $Z(t)=\int_{0}^t g(s)ds$. Now, $g$ can
be seen as the natural parametrisation of $B$, which leads
$g(s)=\exp(iG(s))$ for some function $G:[0,1]\to [0,2\pi)$, non
decreasing. Hence $G$ has a right continuous modification $\tilde G$
which also satisfies $Z(t)=\int_{0}^t e^{i\tilde{G}(s)}ds$. The
function $\tilde{G}$ is
the inverse of a CDF $F_\nu$ for some $\nu$ in~$\MT^0$.\\

\Proofof{Theorem \ref{th:mes-conv} 2)} 
Consider first the continuity of $\BB$. For any $t\in[0,2\pi)$ and any pair of distributions $(\mu,\nu)$, since $x\to \exp(ix)$ is 1-Lipschitz, 
\be
|Z_\mu(t)-Z_\nu(t)|&=&\l|\int_{0}^t \exp(iF^{-1}_\mu(u))-\exp(iF^{-1}_\nu(u)) du \r|\\
&\leq &\int_{0}^t d_\T(F^{-1}_\mu(u),F^{-1}_\nu(u)) du, \ee
where $d_\T$ is the distance in $\T$, defined for $0\leq x\leq y
<2\pi$ by $d_\T(x,y)=\min\{y-x,2\pi-y+x\}$. This last quantity is then bounded above, uniformly in $t\in[0,1]$ by $`E(d_\T(X_\mu,X_\nu))$, for
\[X_\mu:= F^{-1}_\mu(U),~~ X_\nu:=F^{-1}_\nu(U),\] where $U\sim\uniform[0,2\pi]$. Now, $`E(d_\T(X_\mu,X_\nu))$ is a Wasserstein like distance $W_1(\mu,\nu)$ between the distributions $\mu$ and $\nu$
in $\T$ (the standard Wasserstein distance is rather defined between measures on an interval, not on the circle). Now, it is classical that the convergence in
distribution implies the convergence of the Wasserstein distance to~0 (see Dudley \cite{Dud} Section 11.8). 
This property can be easily extended to the present case, considering
that $X_n\dd X$ in $\MT$ iff there exists $\theta\in[0,2\pi]$ (any
point of continuity of $X$ does the job) for which $X_n-\theta \mod
2\pi\dd X-\theta \mod 2\pi$ in the standard sense.

\noindent Reciprocally, let $(B_n,n\geq 0)$ be a sequence of CCS boundaries $B_n$ converging to $\BB_\mu$ for the Hausdorff distance $d_H$. By Theorem~\ref{th:mes-conv}~1), there exists  $\mu_n\in\MT^0$ such that $\BB_{\mu_n}=B_n$. We now establish that $(\mu_n,n\geq 0)$ possesses exactly one accumulation point, equal to $\mu$. Consider a subsequence $F_{\mu_{n_k}}$ such that $F_{\mu_{n_k}}
\overset{D_1}{\longrightarrow} G$, where $G$ is the CDF of a measure $\nu$. Such a subsequence exists since $\MT^0$ is compact (and then sequentially compact, since it is a metric space). Now, for $D_1$ denoting the Skorokhod distance (see e.g. Billingsley \cite{BIL} Chap.3),
$
F_{\mu_{n_k}} \overset{D_1}{\longrightarrow} G \quad \Rightarrow \quad 
F_{\mu_{n_k}}^{-1} \overset{D_1}{\longrightarrow} G^{-1}.
$
According to the first part of this proof, the limit CCS boundary
$\BB_{\nu}$ must be equal to $\BB_\mu$. Since by
Theorem~\ref{th:mes-conv}~1), the CCS characterise the measure,
$\nu\overset{(d)}{=}\mu$.  ~$\Box$ \medskip

\subsection{Fourier decomposition of the CCS curve}
\label{sec:fourier}

Fourier coefficients provide powerful tools to analyse the geometrical
properties of the CCS curves.

\noindent Let $f$ be a function from $[0,2\pi]$ with values in
$`R$. The quantity $\frac{1}{2}a_0+\sum_{k\geq 1}
a_k\cos(ku)+b_k\sin(ku)$ is the standard Fourier series of $f$, where

$$
  \displaystyle \quad a_k=\pi^{-1}\int_0^{2\pi} \cos(ku)f(u)du,
  \quad b_k=\pi^{-1}\int_0^{2\pi} \sin(ku)f(u)du.
$$
For $\mu$ in $\MT$ (or in
${\cal M}[0,2\pi]$), the Fourier coefficients of $\mu$ are defined,
for any $k\geq 0$ by

\begin{equation}\label{eq:co-a-b}
  a_0(\mu) =  \frac{1}{\pi},
\quad
  a_k(\mu)  = \frac{1}{\pi} \mathbb{E}(\cos(kX_\mu)), 
\quad 
  b_k(\mu)  = \frac{1}{\pi} \mathbb{E}(\sin(kX_\mu)). 
\end{equation}
In this setting, the condition $\int_0^{2\pi} e^{iu}dF_\mu(u)=0$
coincides with
\begin{equation}
a_1(\mu)=`E(\cos(X_\mu))=0,\quad b_1(\mu)=`E(\sin(X_\mu))=0.
\end{equation}
The following proposition, whose proof can be found in Wilms \cite[Theorem 1.6 and 1.7]{WilmsFractional}, states that probability measures are characterised by their Fourier coefficients, and establishes a continuity theorem.
\begin{prop}\label{pro:mes-ser}%
1) The function
$$\app{\Coeffs}{\MT}{`R^{\mathbb{N}}\times`R^{\mathbb{N}}}{\mu}{\l((a_k(\mu),k\geq
  0),(b_k(\mu),k\geq1)\r)}$$ 
is injective.\\
2) Let $\mu,\mu_1,\mu_2,\dots$ be a sequence of measures in $\MT$. The two following statements are equivalent: $\mu_n\weak\mu$ and $\Coeffs(\mu_n)$ converges pointwise to $\Coeffs(\mu)$ (meaning that for any $k$, $a_k(\mu_n)\to a_k(\mu)$ and $b_k(\mu_n)\to b_k(\mu)$).
\end{prop}

\begin{exm}\label{exm:reg-gone}
  -- If $\mu\sim \uniform[0,2\pi]$ then
  $a_k(\mu)=b_k(\mu)=0$ for any $k\geq 1$.\\
  -- If $\mu=\sum_{k=0}^{m-1}\frac{1}{m} \delta_{2\pi k/m}$ is the
  uniform distribution on the vertices of a regular $m$-gon (with a
  vertex at position $(0,0)$), then all the $b_k$ are null, 
     $a_0(\mu)=1/\pi$, and
     $a_k(\mu) = \pi^{-1} \1_{k\in m\mathbb{N}^\star}$.
\end{exm}

Of course, deciding whether a given pair $((a_k,k\geq 0),(b_k,k\geq
1))$ corresponds to a pair
$((a_k(\nu),k\geq 0),(b_k(\nu),k\geq 1))$ for some $\nu\in\MT$ is a difficult task: there 
does not exist in the literature any characterisation of Fourier series of non negative
measures. The case of measures having a density with respect to the
Lebesgue measure is discussed in Section~\ref{sec:positive-fourier-series}.

The area of a CCS $\CC_\mu$ has an expression in terms of $\Coeffs(\mu)$. In this section, we consider
a CCS with a smooth $C^1$ boundary that is equal to its Fourier expansion.  
The following formula can be deduced from Hurwitz
\cite[p.372-373]{Hurwitz}, where it is given using a parametrisation
of the boundary of the CCS. In our settings, writing ${\cal
  A}(\mu)$ for the area of $\CC_\mu$, it translates into:
  \begin{equation}\label{eq:area}%
    {\cal A}(\mu) = \frac{1}{4\pi} - \frac{\pi}{2} \sum_{k\geq 2} \frac{a_k^2(\mu)+b^2_k(\mu)}{k^2-1}.
\end{equation}
As did Hurwitz, this equation can be proved from Green's theorem stating that:
\begin{equation}\label{eq:int-are}  {\cal A}(\mu) =\int_0^1 x_\mu(t)\frac{dy_\mu(t)}{dt}dt=-\int_0^1 y_\mu(t)\frac{dx_\mu(t)}{dt}dt.
\end{equation}
As a matter of fact, this formula remains valid for every CCS in
$\Conv(1)$ (cf. Corollary~\ref{cor:general_area_formula}). 
Rewriting \eref{eq:int-are} and using \eref{eq:curve-x-y} gives
\ben
\nonumber  {\cal A}(\mu) &=&\int_0^1 \int_{0}^t \cos(F_\mu^{-1}(u))du \sin(F_\mu^{-1}(t))dt\\
\label{eq:aire1}              &=& `E\l(\cos(X)\sin(X')1_{X\leq X'}\r).
\een
where $X$ and $X'$ are two independent copies of $X_\mu$.
\begin{rem}One can show that \eref{eq:area} implies \eref{eq:aire1} by
  noticing that $`E(\cos(kX))^2+`E(\sin(kX))^2=`E(\cos(k(X-X'))$ and
  using the general equality $ \sum_{k\geq 2} \frac{\cos(kx)}{k^2-1}=
  \frac{\cos(x)}{4}-\frac{(\pi-(x
    \bmod{2\pi}))}2\sin(x)+\frac12$. 
  Notice that Hurwitz~\cite{HurwitzIsoperimetres} deduced the
  isoperimetric inequality from \eref{eq:aire1} with a proof which
  only requires an equivalent of Wirtinger's inequality.
\end{rem}

\subsection{Convergence of discrete CCS and an application to
  statistics}
\label{sec:conver-conv} 

Consider $X_1,\dots,X_n$ \iid having distribution $\mu$ with support in $[0,2\pi)$.  The empirical CDF associated with this sample is defined by $F_n(x)=n^{-1}\#\{i: X_i \leq x\}$. The law of large number ensures that $F_n\to F_\mu$ pointwise in probability, and  $(n^{1/2} |F_n(x)-F_\mu(x)|,x\in[0,2\pi])$ converges in distribution in  $D[0,2\pi]$, the set of càdlàg function equipped with the Skorokhod topology, to $({\sf b}(F_\mu(x)),x \in[0,2\pi])$ where ${\sf b}$ is a standard Brownian bridge (see Billingsley \cite[Theorem 14.3]{BIL}).

Now assume that the $X_i$ take their values in $\T$, and let $\hat
X_1,\dots, \hat X_n$ be the sequence $X_1,\dots,X_n$ sorted in
increasing order (with the natural order on $[0,2\pi)$). Consider
the function $Z_n:[0,1]\to \mathbb{C}$ defined by $Z_n(0)=0$,
\[Z_n(k/n)=\frac{1}{n}\sum_{j=1}^k \exp(i \hat X_j),~~ \textrm{ for }k
\in\{1,\dots,n\},\]
and extended by linear interpolation between the points $(k/n,k\in
\{0,\dots,n\})$. Also define the empirical curve $B_n$ associated with
the distribution $\mu$, as $B_n := \{Z_n(t),t\in[0,1]\}$. The
curve~$B_n$ belongs to $\BConv(1)$ if and only if $\sum_{j=1}^n
e^{iX_j}=0$; otherwise, since the steps are sorted, $B_n$ is either
simple or may contain at most 1 self-intersection point, that is a
pair $t_1<t_2$ such that $Z_n(t_1)=Z_n(t_2)$. For $\theta\in[0,2\pi)$,
let $N_n(\theta) = \#\{i, X_i \leq \theta\}$ be the number of
variables smaller than~$\theta$.  The set of extremal points of $B_n$
is
\beq
\Ext(B_n)=\l\{Z_n(N_n(\theta)/n),\theta\in[0,2\pi]\r\}.
\eq
Set for any $\theta\in[0,2\pi)$, 
\[W_n(\theta):=\sqrt{n}\l[Z_n(N_n(\theta)/n)-Z_\mu(F_\mu(\theta))\r].\]
This process measures the difference between $Z_n$ and its limit.\par
\begin{figure}[tbp]
  \begin{center}
    \includegraphics[width=15.5cm]{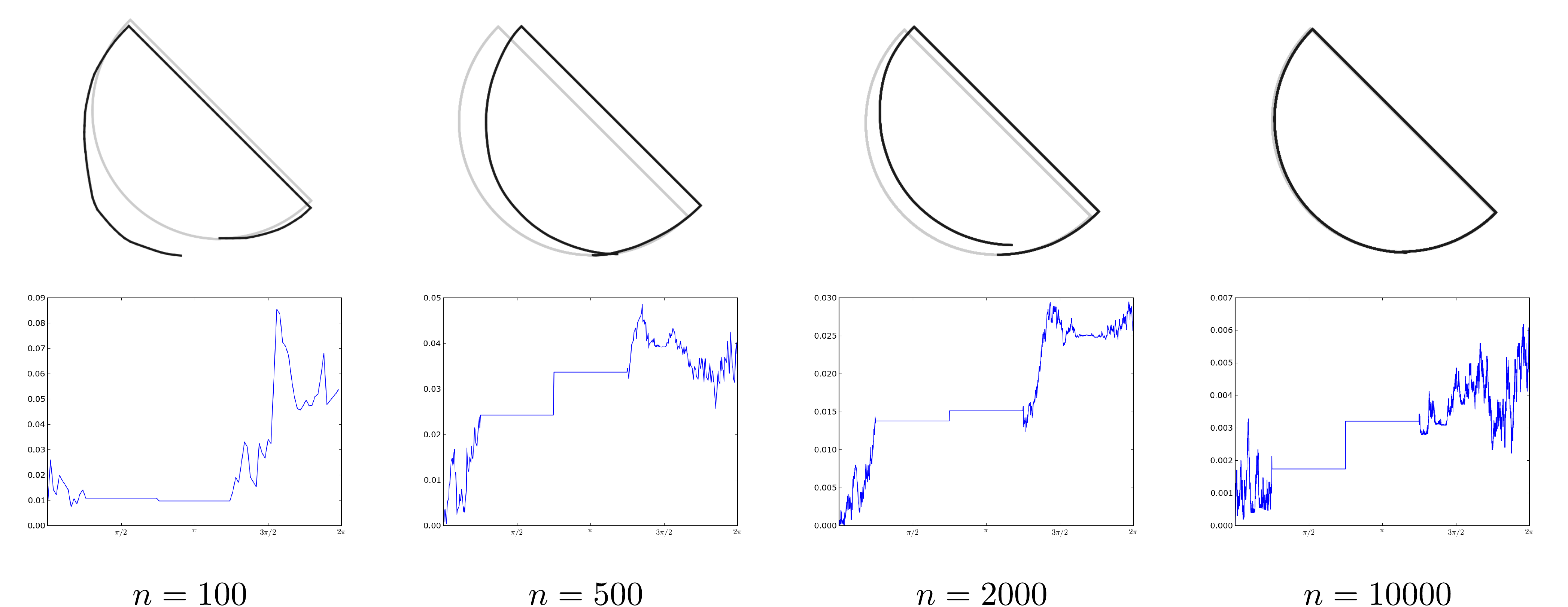}
    ~\vskip -1.2cm~
  \end{center}
  \caption{Convergence towards the half-circle. The first row of
    figures describes the discrete CCS of size $n$ (in black) compared
    to the limit CCS (in grey). The second row displays the distance
    between the discrete CCS and its limit ($\theta \rightarrow
    |W_n(\theta)|$).}
  \label{fig:conv22}
\end{figure}

Denote by $\pi_1(z)=\Re(z)$, $\pi_2(z)=\Im(z)$ and $\pi(z)=(\pi_1(z),\pi_2(z))$.
\begin{theo} \label{thm:cv-emp}
1) The following convergence 
\beq 
\pi \l(W_n(\theta),\theta\in[0,2\pi]\r)\dd (G_\theta,\theta\in[0,2\pi])
\eq
holds in $(D[0,2\pi],`R^2)$, where $G$ is a centred Gaussian process whose finite dimensional distributions are given in Section \ref {proof:thm:cv-emp}, in Formula \eref{eq:FDD}.\\
2) For any $n\geq 1$, $d_H(B_n,\BB_\mu)=\max_\theta |Z_n(N_n(\theta)/n)-Z_\mu(F_\mu(\theta))|$, and then $\sqrt{n}d_H(B_n,\BB_\mu)$ converges in distribution to $\max_\theta |G_\theta|$. 
\end{theo}
See illustration in Figure \ref{fig:conv22}.
The following Corollary -- which gives the asymptotic shape for our random polygons -- is a direct consequence of Theorem \ref{thm:cv-emp}.

\begin{cor}\label{cor:cv-emp}
 If $\mu\in \MT^0$ then:
\begin{itemize*}
\item[1)] The following convergence holds in distribution in $D[0,2\pi]$:
\beq
(Z_n(N_n(\theta)/n),\theta\in[0,2\pi])\dd (Z_\mu(F_\mu(\theta)),\theta\in[0,2\pi]).
\eq
\item[2)] $d_H(B_n,\BB_\mu)\to 0$ in probability.
\end{itemize*}
\end{cor}

\begin{rem}
A direct proof of Corollary~\ref{cor:cv-emp} that ignores Theorem \ref{thm:cv-emp} is as follows: first, the convergence of the finite dimensional distributions (FDD)
corresponding to $1)$ holds as a consequence of the law of large numbers. Then, for an $`e>0$, choose  $k$ and the points $(\theta_1,\dots,\theta_k)$ such that the union of the segments $B_{`e}:=\cup_{i=0..k-1}[Z_\mu(F_\mu(\theta_i)),Z_\mu(F_\mu(\theta_{i+1}))]$ has a length larger than $1-`e$. From there, 2) follows since for $n$ large enough, $|Z_n(N_n(\theta_i)/n)-Z_\mu(F_\mu(\theta_i))|$ goes to 0 in probability for any $i\leq k$. This implies that the union of the segments 
$B'_n=\cup_i[Z_n(N_n(\theta_i)/n),Z_n(N_n(\theta_{i+1})/n)]$ has total length larger than $1-2`e$ for $n$ large enough, with probability going to 1. Since $B_n$ has length 1, for those same $n$, $d_H(B_n,B'_n)\leq 2`e$.
\end{rem}

The proof of Theorem \ref{thm:cv-emp} is postponed to the appendix.

\section{Operations on measures and on CCS}

Mixture and convolution are natural operations on $\MT^0$:
\begin{itemize*}
\item[1)] \textit{Mixture}: if $\mu,\nu \in \MT^0$  then for any $\lambda\in[0,1]$,
$\lambda \mu+(1-\lambda)\nu\in
  \MT^0$.
\item[2)] \textit{Convolution}: if $\mu,\nu \in \MT^0$ then $\mu \CT \nu \in
  \MT^0$, where $(\CT)$ denotes the convolution in $\MT$. This
  conclusion holds even if only $\mu$ is in $\MT^0$.
\end{itemize*}

\noindent Then the maps $\BB$ and $\CC$ transport
these operations on $\Conv(1)$:
\begin{defi}\label{def:convo}
  Let $\CC_\mu$ and $\CC_\nu$ be two CCS in $\Conv(1)$ and  $\lambda\in[0,1]$.
  \begin{itemize*}
  \item[1)] We call mixture of $\CC_\mu$ and of $\CC_\nu$ with weights
    $(\lambda,1-\lambda)$, the CCS $\CC_{\lambda\mu+(1-\lambda)\nu}$.
  \item[2)] We call convolution of $\CC_\mu$ and $\CC_\nu$, the CCS
    $\CC_\mu \star\CC_\nu:=\CC_{\mu\CT\nu}$.
  \end{itemize*}
\end{defi}
In this section we provide some facts which seem to be unknown: a
mixture is sent by $\CC$ on a Minkowski sum (Proposition
\ref{prop:mixture_minkowski}) and the Minkowski symmetrisation can
also be expressed in terms of mixtures (Theorem \ref{thm:pro-sym}).
The convolution of CCS acts somehow on the radius of curvature and
seems to be a new operation, leading to a notion of symmetrisation by
convolution that we introduce in section \ref{sec:Conv}.

\begin{figure}[htbp]
  \begin{center}
    \begin{tabular}{cc}
      ~\includegraphics[width=6.5cm]{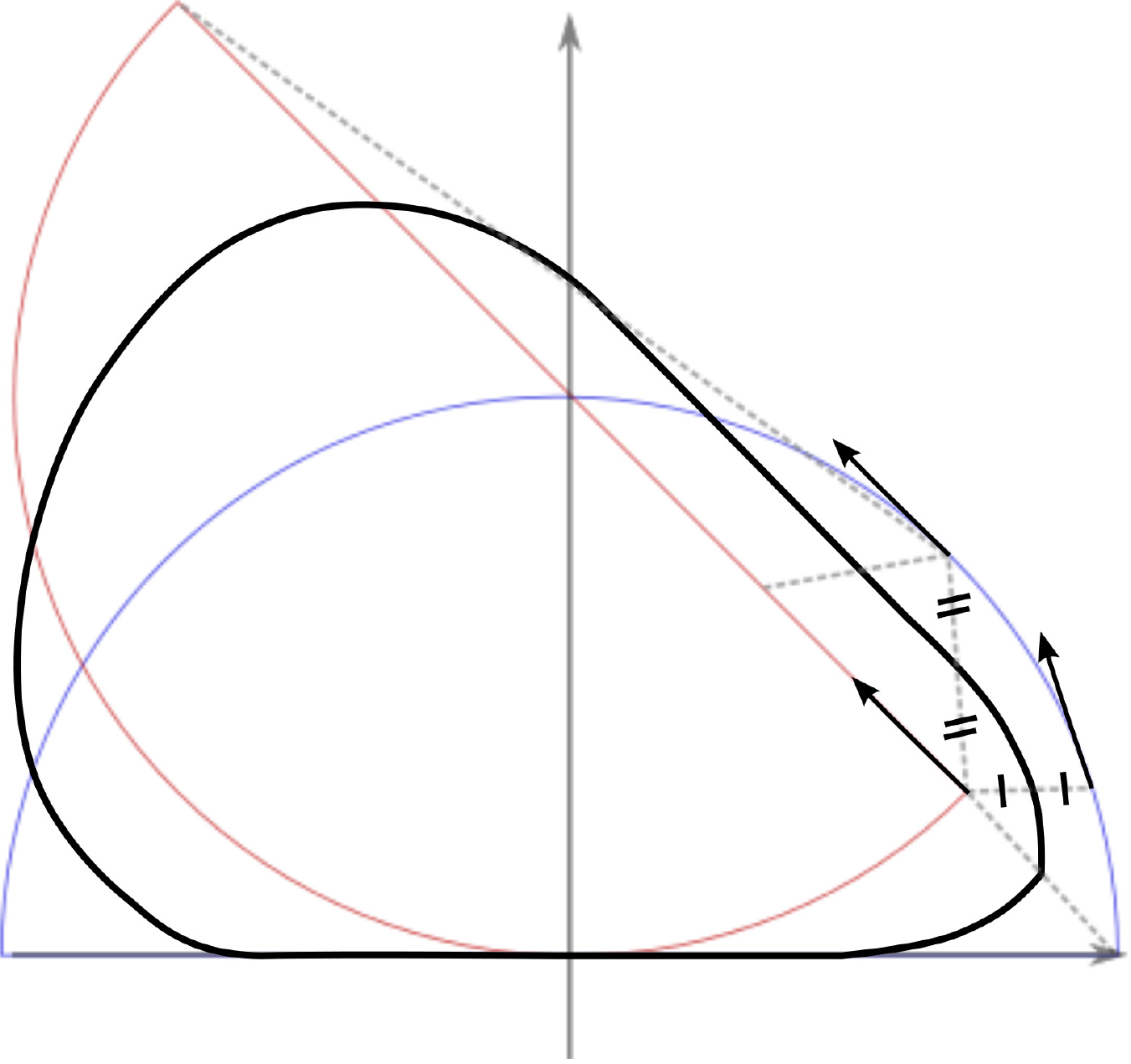}~ & 
      ~\includegraphics[width=6.5cm]{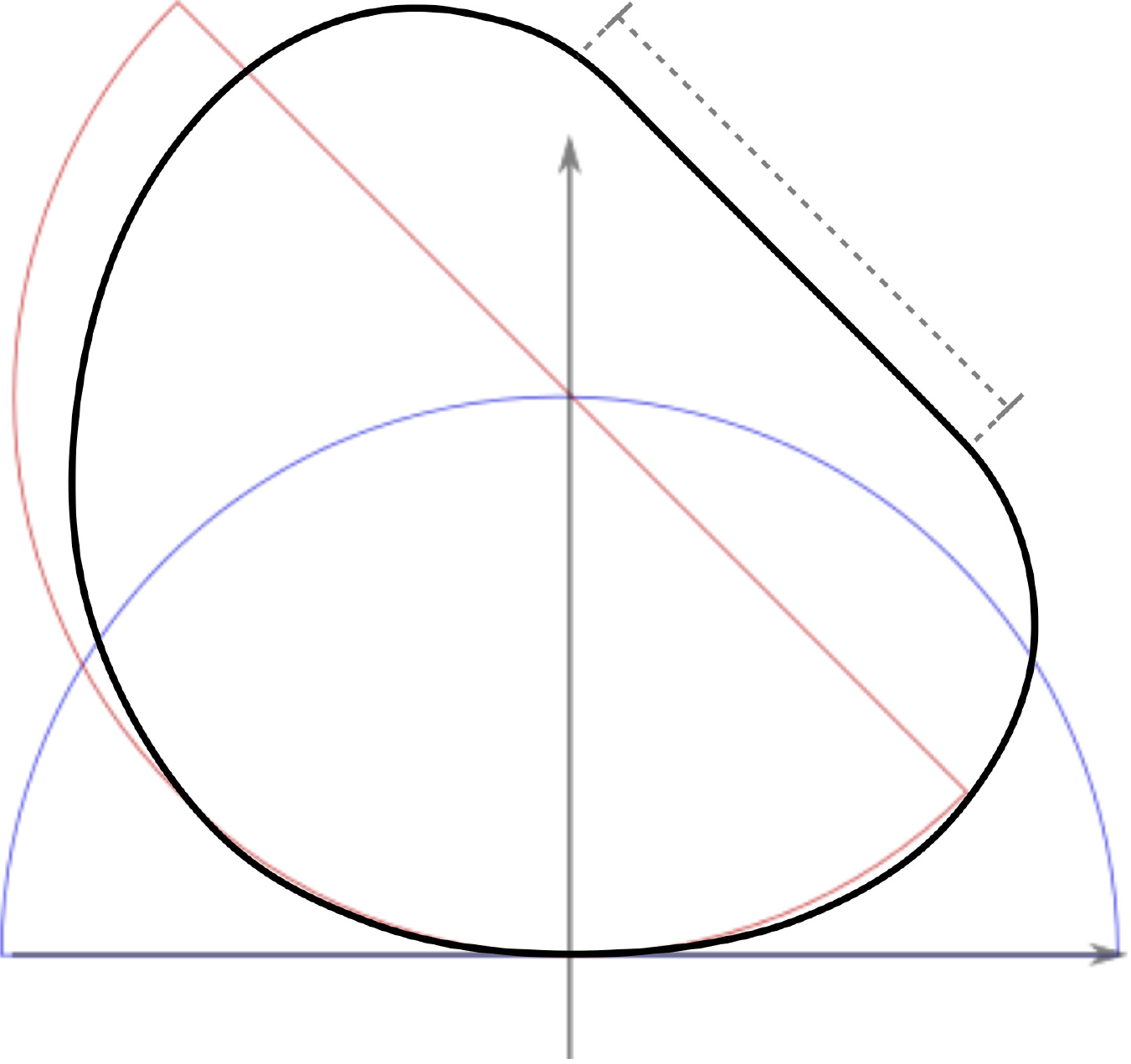}~ \\
      {\small (a)} & {\small (b)}
    \end{tabular}
    ~\vskip -1cm~
  \end{center}
  \caption{\label{fig:conv}Construction of the (a) mixture and (b)
    convolution of two half-circles. Notice that every point of the
    mixture is the barycentre of two points of the original
    half-circles, and that the CCS obtained by convolution
    possesses a linear segment whose angle corresponds to the sum of
    the angles of the segments in the original half-circles.}
\end{figure}

\subsection{Mixtures of CCS / Minkowski sum}
\label{sec:Mixt} 

Let $A$ and $B$ be two subsets of $`R^2$. The Minkowski sum of $A$ and
$B$ is the set $ A+B=\{a+b~: a\in A,b\in B\}$.  Further, for any
$\lambda$, write $\lambda A=\{\lambda a~:a\in A\}$. We have:
\begin{prop}\label{prop:mixture_minkowski}
Let $\nu,\mu\in \MT^0$, $\lambda\in[0,1]$. Then 
\[
\CC_{\lambda \mu+(1-\lambda)\nu}=\lambda \CC_{\mu}+(1-\lambda)\CC_\nu
\]
which means that the mixture of CCS and the Minkowski sum are the
same, and that the CCS of a mixture corresponds to the mixture of the
CCS.
\end{prop}  

This proposition (see Figure \ref{fig:conv}) implies that the boundaries
verify:
\[
\BB_{\lambda \mu+(1-\lambda)\nu}=\partial(\textrm{convex hull}(\lambda \BB_{\mu}+(1-\lambda)\BB_\nu))
\]

\Proof %  
We first give a proof when  $\mu$ and $\nu$ have densities.
Recall the characterisation given in Lemma
\ref{lem:ext-points}. Write
\ben
Z_{\lambda\mu+(1-\lambda)\nu}(F_{\lambda\mu+(1-\lambda)\nu}(t)) 
& = & 
\lambda \int_0^t \exp(it) d\mu(t)  +(1-\lambda) \int_0^t \exp(it)
d\nu(t) \notag \\
& = & \lambda Z_\mu(F_\mu(t))+(1-\lambda)Z_\nu(F_\nu(t)).\label{eq:mix}
\een 

The extremal points of $\CC_{\lambda\mu+(1-\lambda)\nu}$ are then
obtained as particular barycentres of extremal points of $\CC_\mu$ and
$\CC_\nu$. 
When both $\mu$ and $\nu$ have a density, this implies that the point
in $\BB_{\lambda\mu+(1-\lambda)\nu}$ where the tangent has direction
$\theta$ is obtained as the barycentre of the corresponding points in
$\BB_\mu$ and $\BB_\nu$. This implies that
$\CC_{\lambda\mu+(1-\lambda)\nu} \subset \lambda
\CC_{\mu}+(1-\lambda)\CC_\nu$.\par

We establish the other inclusion by
using the fact that CCS are characterised by their supporting
half-planes: for every $t\in[0,2\pi]$, let $D_\mu(t)$ be the line
passing through $Z_\mu(F_\mu(t))$ making an angle $t$ with the
$x$-axis.  The line $D_\mu(t)$ defines a supporting half-plane
$H_\mu(t)$ for $\CC_\mu$.  Since $\CC_\mu$ is a CCS, this half-plane
is minimal for the inclusion with regard to the property of making
an angle $t$ with the $x$-axis. Considering that the points in
\eref{eq:mix} all belong to their associated half-plane, these
half-planes verify:
$$H_{\lambda\mu+(1-\lambda)\nu}(t) = \lambda H_\mu(t) + (1-\lambda) H_\nu(t).$$
Now, the left-hand side represents a supporting half-plane for
$\CC_{\lambda\mu+(1-\lambda)\nu}$ and the right-hand side another
supporting half-plane for $\lambda \CC_\mu + (1-\lambda) \CC_\nu$. We
deduce that the CCS they enclose are equal.\par
When $\mu$ or $\nu$ have no densities, take a sequence $(\mu_n,\nu_n)$ of measures having densities and which converges weakly to $(\mu,\nu)$; we then obtain $\CC_{\lambda \mu_n+(1-\lambda)\nu_n}=\lambda \CC_{\mu_n}+(1-\lambda)\CC_{\nu_n}$ and conclude by Theorem \ref{th:mes-conv}.  ~$\Box$ \medskip

\noindent Hence the CCS $\CC_{\lambda\mu+(1-\lambda)\nu}$ has a perimeter equal to~1, as all CCS of $\Conv(1)$. This implies that the perimeter of the Minkowski sum $\lambda
\CC_{\mu}+(1-\lambda)\CC_\nu$ is 1 (well known fact, obtained here without geometric arguments).
\begin{rem}For $\mu$ and $\nu$ in $\MT^0$ and $\lambda\in[0,1]$, we have 
\beq\label{eq:Mixture1}
{\mathcal A}(\lambda \mu+(1-\lambda)\nu)^{1/2}\geq \lambda {\mathcal A(\mu)}^{1/2}+(1-\lambda){\mathcal A(\nu)}^{1/2}.
\eq This is the so-called Brunn-Minkowski inequality; it implies that ${\mathcal A}(\lambda \mu+(1-\lambda)\nu)\geq \min \{{\mathcal A}(\mu),{\mathcal A}(\nu)\}$. It can be proved using Hurwitz formula \eref{eq:area} and the Cauchy-Schwarz inequality.
\end{rem} 

\subsubsection{Minkowski symmetrisation and measure symmetrisation}

Let $K$ be a CCS of $`R^2$ and $u\in `R^2$, $|u|=1$.  We denote by
$\pi_u \in O(2)$ the reflection with respect to the straight line
passing through the origin and orthogonal to $u$, i.e. $\pi_u(x) = x -
2\langle x, u\rangle u$.  The \textit{Minkowski} (or Blaschke)
\textit{symmetrisation} of $K$ is the CCS
$S_u(K)=\frac12(\pi_uK+K)$. The same operation can be defined over
$\mathbb{C}$: for $u= e^{i\theta}$, the Minkowski symmetrisation of
$K$ with respect to direction $\theta$ is the map $(K,\theta)\mapsto
\frac{e^{i\theta}}{2} (\overline{e^{-i\theta}K}+e^{-i\theta} K)$,
where $\bar{z}$ is the complex conjugate of $z$.

Now, let $\theta\in[0,2\pi]$, $\mu\in\MT^0$, and set $\mu(\theta)$ be
the distribution of $X_\mu+\theta\mod 2\pi$.  Since
$`E(\exp(i(X_\mu+\theta)))=e^{i\theta}`E(\exp(iX_\mu))$, $\mu(\theta)$
is in $\MT^0$. The CCS $\CC_{\mu(\theta)}$ can be obtained from
$\CC_{\mu}$ by a rotation (of angle $-\theta$) followed by a
translation.

For any $\nu\in\MT^0$, set $\overleftarrow{\nu}=\nu(2\pi-.)$.  The
symmetrisation of $\nu$ with respect to direction $\theta$ is the
measure $S(\nu(\theta))$ defined by 
\beq\label{eq:Snu}
S(\nu(\theta))=\frac{1}{2}(\nu(\theta) +\overleftarrow{\nu(\theta)}).
\eq 
Further the symmetrisation by mixture of $\CC_\nu$ with respect to
direction $\theta$ is defined to be $\CC_{S(\nu(\theta))}$.

A direct consequence of Proposition \ref{prop:mixture_minkowski} is
the following:
\begin{prop}\label{pro:mix=Min}%
  The symmetrisation by mixture with respect to direction
  $\theta$ coincides with the Minkowski symmetrisation with
  respect to $u=e^{i\theta}$.
\end{prop}
Again Theorem \ref{thm:cv-emp} provides a new point of view on this
symmetrisation. Starting from a set of angles
$\theta_1,\dots,\theta_k$ and an initial measure $\nu \in \MT^0$,
construct the sequence of measures $\nu_k$ defined by $\nu_0=\nu$ and
$\nu_{k+1} = S(\nu_k(\theta_k))$. This sequence consists in
alternating rotations and symmetrisations of the initial measure $\nu$. 

\begin{theo}\label{thm:pro-sym}
For any $\theta\in[0,2\pi]$, any $\nu\in\MT^0$, the following
properties hold:
\begin{itemize*}
\item[1)] the CCS $\CC_{S(\nu(\theta))}$ has the same perimeter as
  $\CC_{\nu}$ (that is 1),
\item[2)] the area does not decrease: $ {\cal A}(S(\nu(\theta)))\geq
  {\cal A}(\nu)$,
\item[3)] for any $k\geq 0$, there exists $\theta_1,\dots,\theta_k\in[0,2\pi]$ such that
  $$d_H(\CC_{\nu_k},\Circle(i/(2\pi),1/(2\pi)))\leq 2^{-k}\pi,$$
where $\Circle(z,r)$ is the circle with centre $z$ and radius $r$,
\item[4)] among all CCS with perimeter 1, the circle has the largest area.
\end{itemize*}
\end{theo}
Properties 1), 2), 4) are classical; we provide direct probabilistic
proofs below. Statement 3) which gives a bound on the speed of
convergence to the ball for well chosen directions of symmetrisation,
is known in $`R^n$ (see Klartag \cite[Theorem 1.3]{Klartag}), but the
proof we provide here in $`R^2$ is much simpler.

\Proof First, 4) is clearly a consequence of the three first points
(to be honest, our proof uses \eref{eq:Mixture1}, which implies
directly the isoperimetric inequality). The first item follows from
the fact that if $S(\nu(\theta))\in \MT^{0}$, then
$\BB_{S(\nu(\theta))}\in\BConv(1)$. And \eref{eq:Mixture1} implies 2)
since ${\cal A}(\nu)={\cal A}(\nu(\theta))={\cal
  A}(\overleftarrow{\nu(\theta)})$.

Let us prove 3). If $L=[X_1,\dots,X_l]$ for some $l\geq 1$, a list of
\rv with distribution $\nu_1,\dots,\nu_l$, we say that $\nu$ is the
\textit{equi-mixture} of $L$ if $\nu=\frac{1}{l}(\nu_1+\dots+\nu_l)$.

Take $X\sim \nu$. $\nu_1 := S(\nu(\theta_1))$ is the equi-mixture of
$[X+\theta_1 \mod 2\pi,-X-\theta_1 \mod 2\pi]$.  Therefore using that
$(a \mod 2\pi) +b \mod 2\pi =(a+b)\mod 2\pi$, $S_{\nu_2}$ is the
equi-mixture of $[X+\theta_1\pm\theta_2 \mod 2\pi, -X-\theta_1
\pm\theta_2 \mod 2\pi]$.  Iterating this, one observes that
$S_{\nu_k}$ is the equi-mixture of $[X+\theta_1\pm\theta_2 \pm\dots\pm
\theta_{k}\mod 2\pi, -X-\theta_1 \pm\theta_2 \pm\dots\pm
\theta_{k}\mod 2\pi]$.  If $\theta_k=(2\pi)/2^{k-1}$ then $S_{\nu_k}$
is the equi-mixture of $\mu_1$ and $\mu_2$, where $\mu_1$ and $\mu_2$
are the respective equi-mixture of $[X+\theta_1\pm\theta_2 \pm\dots\pm
\theta_{k}\mod 2\pi]$ and of $[-X-\theta_1 \pm\theta_2 \pm\dots\pm
\theta_{k}\mod 2\pi]$.

Now, both $\mu_1$ and $\mu_2$ converge to $\uniform[0,2\pi]$: to check
this, consider the sequence of intervals $I_n=[2\pi n2^{-{k-1}},2\pi
(n+1)2^{-{k-1}})$ , for $0\leq n\leq 2^{k-1}-1$.  For $j\in\{1,2\}$,
$\mu_j(I_n)=1/2^{k-1}$ for any~$n$. Indeed, $\mu_1$ (resp. $\mu_2$) is
the equi-mixture of all measures obtained from the distribution of $X$
(resp. $-X$) by dyadic translation of depth $k$, then since all
intervals $I_n$ have depth $k$, they have the same weight.  Hence
$F_{\mu_1}(2\pi n2^{-k+1} )=n2^{-k+1}$ for any $n$. Therefore, since
$F_{\mu_1}$ is increasing, we have that $\|F_{\mu_j}-F\|_{\infty}\leq
2^{-k+1}$, for $F_\upsilon(x)=x/(2\pi)$, the CDF of
$\uniform[0,2\pi]$, which gives $\|F_{\nu_k}-F_\upsilon\|_{\infty}\leq
2^{-k+1}$. 
Further, the right inverses $F^{-1}_{\nu_k}$ and
$F^{-1}_\upsilon$ are close:
\[\|F^{-1}_{\nu_k}-F^{-1}_\upsilon\|_{\infty}\leq 2^{-k+1} 2\pi.\]
Thanks to \eref{eq:conv-measure2}, 
\be
|Z_{\nu_k}(t)-Z_{\upsilon}(t)|&\leq & \int_0^t\l| \exp(iF_{\nu_k}^{-1}(u))-\exp(iF_{\upsilon}^{-1}(u))\r| du\\
&\leq & \int_0^t\l| F_{\nu_k}^{-1}(u) -F_{\upsilon}^{-1}(u))\r| du 
\ee
and therefore $\|Z_{\nu_k}(t)-Z_{\upsilon}(t)\|_{\infty}\leq 2^{-k}\pi$. ~$\Box$

\subsection{Convolution of measures / Convolution of CCS}
\label{sec:Conv}

In fact, $\BB_{\mu\CT \nu}$ is obtained as a kind of convolution of
$\BB_\mu$ and $\BB_\nu$. As seen earlier if $\mu$ has a density
$f_\mu$ then $f_\mu(\theta)$ represents the radius of curvature of
$\BB_\mu$ at time $F_\mu(\theta)$. Therefore the radius of curvature
$R_\theta$ of $\BB_{\mu\CT \nu}$ at time $F_{\mu\CT \nu}(\theta)$ is
the convolution of the radii of curvature of~$\BB_\mu$ and~$\BB_\nu$
as follows:
\[ R_\theta = \int_0^{2\pi} f_\mu(x) f_\nu((\theta-x)\mod 2\pi) dx. \]
\begin{theo}\label{thm:Mixture2}%
Let $\mu$ and $\nu$ in $\MT^0$. The convolution does not decrease the area
\[{\cal A}\l(\mu\CT\nu\r)\geq \max\{{\cal A}(\mu),{\cal A}(\nu)\}.\]
Since $\uniform[0,2\pi]$ is an absorbing point for $\CT$, and $\CC_u$
is the circle of perimeter $1$, this implies the isoperimetric
inequality: $\mathcal{A}(\uniform[0,2\pi]) \geq \cal{A}(\nu)$,
$\forall \nu \in \MT^0$.
\end{theo} 
\Proof Consider $X$ and $Y$ two independent \rv such that $X\sim \mu$,
$Y\sim \nu$. Let $\eta=\mu\CT \nu$. By expansion of $\cos(n(X+Y))$ and
$\sin(n(X+Y))$ we get
\be
 a_n(\eta)&=& a_n(\mu)a_n(\nu)- b_n(\mu)b_n(\nu)\\
 b_n(\eta)&=& b_n(\mu)a_n(\nu)+ a_n(\mu)b_n(\nu).
\ee
Since $\cos(kX)$ and $\sin(kX)$ have non-negative variances, 
\[a_n^2(\mu)+b_n^2(\mu) =`E(\cos(nX))^2+`E(\sin(nX))^2\leq `E(\cos^2(nX)+\sin^2(nX))=1.\]
Hence,\be
a_n^2(\eta)+b_n^2(\eta)&=&(a_n^2(\mu)+b_n^2(\mu))(a_n^2(\nu)+b_n^2(\nu))\\
                       &\leq&  \min\{a_n^2(\mu)+b_n^2(\mu),a_n^2(\nu)+b_n^2(\nu)\},
\ee
The conclusion follows from \eref{eq:area}.~$\Box$ 

\begin{cor} \label{cor:general_area_formula}%
  Let $\mu \in \MT^0$. Then the formula~(\ref{eq:area}) for ${\cal A}(\mu)$ holds.
\end{cor}
\noindent\Proof
Formula~(\ref{eq:area}) is valid when $\mu$ admits a $\mathcal{C}^1$ density. Just assume that $`E(e^{iX_\mu})=0$.
Let $N$ be a Gaussian centred \rv with variance 1, and let $N_k=N/\sqrt{k} \mod 2\pi$ for $k\geq 1$, and $\mu_k = \mu * N_k$. Clearly $\mu_k\in\MT^0$, and $\mu_k\weak \mu$ which implies ${\cal A}(\mu_k)\to {\cal A}(\mu)$. Now,
\[ \forall n \in \Rel, \qquad  
`E(e^{i n N_k}) = `E(e^{i n (N/\sqrt{k} \mod 2\pi)}) =`E(e^{i n N/\sqrt{k}})=e^{-\frac{n^2}{2k}}.\]
Then the Fourier coefficients of $N_k$ verify 
$a_n = e^{-\frac{n^2}{2k}}$ and $b_n = 0$. Since $\mu_k$ admits a $\mathcal{C}^\infty$ density function, and as a corollary of the proof of Theorem~\ref{thm:Mixture2}: 
$$ {\cal A}(\mu_k) = \frac{1}{4\pi} - \frac{\pi}{2} \sum_{n\geq 2}
\frac{(a_n^2\left(\mu) + b_n(\mu)^2\right)  e^{-\frac{1}{2k}n^2}}{n^2-1}. $$
As a consequence of Lebesgue's dominated convergence theorem, ${\cal
  A}(\mu_k)$ converges to the right hand side of \eref{eq:area}.
~$\Box$

\begin{defi} A measure $\nu$ in $\MT$ is said to be $c$-stable (for some $c>0$) if for $X_\nu$ and $X'_\nu$ two independent \rv under $\nu$,
\ben
\label{eq:c-stable} X_\nu+X'_\nu \mod 2\pi \sur{=}{(d)} c X_\nu
\mod 2\pi.
\een \end{defi} 
This qualification of ``stable'' comes from the standard notion of probability theory where the same question is studied without the $\bmod{\ 2\pi}$ operation (see Feller \cite[Section VI]{FEL2}).\par
The following Proposition due to Lévy \cite[p.11]{LP} identifies  the set of $1$-stable distributions.
\begin{prop}The only 1-stable measures are $\uniform[0,2\pi]$, the Dirac measure at 0, and the family, indexed by $m\geq 1$, of uniform measures on $\{k 2\pi/m,  k=0,\dots,m-1\}$.
\end{prop}

We say that a distribution $\nu$ is in the $2\pi$-domain of attraction of a distribution $\mu$, and write $\nu\in \DA(\mu)$, if for a family $(X_i,i\geq 1)$ of \iid \rv under $\nu$, there exists $\theta\in[0,2\pi]$ such that 
\[\sum_{i=1}^n(X_i-\theta) \mod 2\pi \dd X_\mu.\] 
We let $\DA=\{\mu~: \DA(\mu)\neq \emptyset\}$ be the set of measures $\mu$ whose domains of attraction are not empty.
\begin{prop}\label{pro:domain-att} 
1) The set $\DA$ coincides with the set of 1-stable distributions. \\
2) For any $\nu\in\MT^0$, there exists $\theta\in[0,2\pi]$ and a unique 1-stable measure $\mu$ s.t. $\nu \in \DA(\mu)$.
\end{prop} 
\Proof  1) If $\nu$ is a 1-stable distribution, and if $(X_i,i\geq 1)$ are \iid and taken under $\nu$, then it is easily seen that 
 $X_1+\cdots + X_n \mod 2\pi \sur{=}{(d)} X_1$. Therefore, every 1-stable distribution is in $\DA$. \\
Conversely, assume that $(X_i,i\geq 1)$  are \iid, distributed according to $\nu$, and that $\sum_{i=1}^n(X_i-\theta)\mod
2\pi\dd\mu$. Splitting the sum on the left-hand side into two parts,  $\mu$ appears to be  solution of $\mu=\mu\CT\mu$, and then $\mu$ is 1-stable.\\
2) Take $(X_i,i\geq 1)$ \iid \rv under $\nu$, $\theta\in[0,2\pi]$, and
compute the limit of the $k$-th Fourier coefficient, for $k\geq 1$, of
$\sum_{j=1}^n(X_j-\theta)$,
$$`E(e^{ik\sum_{j=1}^n (X_j - \theta)}) =
`E(e^{ik(X_1-\theta)})^n.$$ 
This coefficient either converges to $0$ or is of modulus $1$ (which
implies $X=\theta/k [2\pi/k]$ \AS). In either case, the limit is a
1-stable distribution. More precisely, let $k$ be the smallest Fourier
coefficient of the limit of modulus $1$. If $k=+\infty$, the limit is
the uniform distribution on $[0,2\pi]$, otherwise it is the uniform
distribution on $\{\frac{2j\pi}{k},j \in [0,k-1]\}$.  (see also Wilms
\cite[Thm. 2.1 and Thm. 2.4]{WilmsFractional}). ~$\Box$

\subsection{Symmetrisation of CCS by convolution} \label{sec:symm}

Let $\nu\in\MT^0$ and $\overleftarrow{\nu}=\nu(2\pi-.)$. The distribution 
\beq\label{eq:sym-conv} %
S_C(\nu):=\nu \CT \overleftarrow{\nu}
\eq is clearly symmetric. We call it the \textit{symmetrisation by convolution} of $\nu$.\footnote{Notice that in the definition of the symmetrisation, replacing $2\pi$ by some other $\theta$ (in $\overleftarrow{\nu}$) affects $S_C(\nu)$ by a simple rotation in $\T$.} 

Denote by $\nu_1=S_C(\nu)$, $\nu_2=S_C(\nu_1)$, ... Let $X_n$ be a \rv under $\nu_n$.
\begin{prop}\label{prop:sym-conv} Let $\nu\in \MT^0$, and let $\mu$ be the unique measure such that $S_C(\nu)$ belongs to
 $\DA(\mu)$. For $\theta=\pi$ or $\theta=0$ we have 
\[X_n-n\theta\mod 2\pi\dd \mu.\]
\end{prop}
\Proof First, $\nu_n$ is the distribution of $\sum_{i=1}^n
(X_i-X'_i)\mod 2\pi$ for some \iid copies $X_i's$ and $X'_i$'s of $X_\nu$. 
The Fourier coefficients of $\nu_n$ can then be computed, and they converge to those of a 1-stable distribution
as in Proposition~\ref{pro:domain-att}, for $\theta\in\{0,\pi\}$ since $X_i-X'_i$ is symmetric.$~\Box$

\section{Extensions}

In this section are discussed two natural extensions of our model. In Section \ref{sec:CWNFP} we discuss CCS with an unconstrained perimeter. In Section \ref{sec:mes-C} is investigated the convergence of a trajectory made by \iid increments with values in $\mathbb{C}$ sorted according to their arguments. If $\nu$ is a centred distribution on $\mathbb{C}$, these trajectories converge to a CCS $\CC_{\sK(\nu)}$ for an operator $\sK$ defined below.
 
\subsection{CCS with an unconstrained perimeter}
\label{sec:CWNFP}
The perimeter of the CCS in the construction we gave is 1 because the total mass of all measures in $\MT^0$ 
is 1. Denote by $\overline{\MT}^0$ the set of positive measures $\nu$ with support $\T$ and such that $\nu(\T)<+\infty$.  Formula 
\eref{eq:conv-measure2}, which defines the CCS associated with a probability measure extends to these 
measures, and the CCS perimeter $\Per(\nu)=\nu(\T)$. A lot of statements given before extend naturally to $\overline{\MT}^0$. Most notably 
\begin{prop} For any measures $\nu_1,\nu_2\in \overline{\MT}^0$, any positive numbers $\lambda_1,\lambda_2$ we have:
\ben\label{eq:perip}
\Per\l(\sum_{i=1}^n \lambda_i\nu_i\r)&=&\sum_{i=1}^n \lambda_i\, \Per(\nu_i)\\
\Per\l(\nu_1\star \nu_2\r)&=&\Per(\nu_1)\, \Per(\nu_2).
\een
The area of $\CC_{\sum_{i=1}^n \lambda_i\nu_i}$ and of $\CC_{\nu_1\star \nu_2}$ are still given by the Fourier coefficients of the measures $\sum_{i=1}^n \lambda_i\nu_i$ and $\nu_1\star \nu_2$, as can be easily checked.
\end{prop}
As said before, \eref{eq:perip} is a well known result.

\subsection{Reordering of random vectors in $\mathbb{C}$} 
\label{sec:mes-C}

The Gauss-Minkowski correspondence can be seen thanks to Corollary \ref{cor:cv-emp} as a consequence of the convergence of polygonal lines corresponding to some reordered random segments. This reordering can be done even if the lengths are not all the same; nevertheless the condition $`E(e^{iX_\mu})=0$ is needed to get a closed convex curve at the limit. In this section we investigate a generalisation of this construction where the sides of the polygons are \rv in $\mathbb{C}$. \par

Let $\mu$ be a distribution with support included in $\mathbb{C}$ with
mean 0, but different from $\delta_0$. Take a sequence
$W:=(W_1,\dots,W_n)$ of \iid \rv with common distribution $\mu$, and
let $\hat{W}:=(\hat{W}_1,\dots,\hat{W}_n)$ the list $W$ sorted
according to the arguments of the $W_i$'s (if several of them have the
same argument but different modulus, then take a uniform random order
among them). For $\theta \in [0,2\pi)$, define $N_n(\theta) := \#\{i,
W_i \leq \theta\}$. Let $S:=(S(k),k=0,\dots,n)$ be the sequence of
partial sums
\beq\label{eq:Sk} S(k):=\sum_{j=1}^k \hat{W}_j,
\eq
piecewise linearly interpolated between integer points, and let
$\bB_n=\{S(t),t\in[0,n]\}$ be the polygonal line corresponding to the
graph of $S$ extended to $[0,n]$. \par

The distribution $\mu$ induces a law $`P_{|W|,\arg(W)}$ for the pair
$(|W|,\arg(W))$, and a law $`P_{\arg(W)}$ for $\arg(W)$; let
$`P_{|W|,x}$ be a version of the distribution of $|W|$ conditioned on
$\arg(W)=x$ (this is defined up to a null set under $`P_{\arg(W)}$;
for the sake of completeness, take $`P_{|W|,x}=\delta_0$ on the
complementary set).  We denote by $m_x$ the mean of $|W|$ under
$`P_{|W|,x}$.

Let $\nu$ be the measure having density $m/`E(|W|)$ with respect to $`P_{\arg(W)}$, that is 
\beq\label{eq:dnu}
d\nu(x)=\frac{m_x}{`E(|W|)}d`P_{\arg(W)}(x).
\eq
The map which sends $\mu$ onto $\nu$ will be denoted $\sK$:
\beq\label{eq:dnu2}\sK(\mu)=\nu.
\eq
Denote by $F^{\arg}$ the CDF of $\arg(W)$, and by $F_\nu$ that of the measure $\nu$.  From now on, let $W_{\theta}$ denote a \rv $W$ under the condition $\{\arg(W)\leq \theta\}$.\par

We here present a theorem stating the aforementioned convergence; we think that it provides an agreeable way to see the phenomenons into play.
\begin{theo}\label{thm:cv-conve} 
Consider the model described in the present section. Assume that $\mu$
is centred ($\neq \delta_0$), 
and let $\nu=\sK(\mu)$. We have\\
1) $d_H(\bB_n/(n`E(|W|)),\BB_{\nu})\as 0$.\\
2) For any $\theta$,  
\beq\label{eq:CVcompl}
\frac{S(N_n(\theta))}{n`E(|W|)}\as  \int_{0}^{\theta}  e^{it} d\nu(t)=Z_\nu(F_\nu(\theta)).
\eq
\end{theo}

\begin{rem} 
  (a) Prosaically, the previous Theorem says that if $\mu$ is a
  centred distribution on~$\mathbb{C}$ the CCS associated with $\mu$
  is $\CC_{K(\mu)}$. \\
  (b) According to \eref{eq:dnu} and Theorem \ref{thm:cv-conve},
  $\BB_{\sK(\nu)}$ is the circle (with radius $1/(2\pi)$) if and only
  if $`P_{\arg}$ admits a density $f_\nu(\cdot)$ with respect to the Lebesgue
  measure, and $\theta\mapsto f_\nu(\theta) m_\nu(\theta)$ is constant. \\
  (c) The ellipse of equation $x^2/c^2 + y^2 = R^2$ with perimeter
  $2\pi R c=1$, is obtained in the case where$$m_\nu(\theta) =
  \frac{1}{2\pi} \frac{c}{\cos(\theta)^2 + c^2\sin(\theta)^2}.$$ This
  can be shown using the following parametrisation: $x(t)=\sin(t)$,
  $y(t)=c(1-\cos(t))$. \\
\end{rem}
\Proofof{Theorem \ref{thm:cv-conve} 2)} 
The cardinality of $N_n(\theta)$ has the binomial
$(n,F^{\arg}(\theta))$ distribution. It satisfies for any $\theta$,
\beq\label{eq:pre}
N_n(\theta)/n \as F_\nu(\theta).
\eq
Conditionally on $N_n(\theta)=m$ the (multi)set $\{\hat{W}_1,\dots,\hat{W}_m\}$ is distributed as a set 
of $m$ \iid copies of $W_\theta$. Therefore by the law of large number,  
\ben\label{eq:deuz}
\frac{S(N_n(\theta))}{n`E(|W|)}&\as&\frac{F^{\arg}(\theta) `E(W_{\theta})}{`E(|W|)}
=\frac{`E(            W1_{\arg(W)\leq \theta})}{`E(|W|)}\\
&=&\frac{`E(|W|e^{i\arg(W)}1_{\arg(W)\leq \theta})}{`E(|W|)}\\
&=& \int_{0}^{\theta}  e^{it} \frac{m_t}{`E(|W|)} d`P_{\arg(W)}(t)=Z_\nu(F_{\nu}(\theta)).
\een
This ends the proof of 2) and shows the \AS simple convergence of the extremal points of the random curve to those of the deterministic limit. \par

\Proofof{Theorem \ref{thm:cv-conve} 1)} Similarly, the length
$L_n(\theta)$ of the curve composed by the segments between the
points $(S(i),0\leq i\leq N_n(\theta))$ satisfies
\ben\label{eq:long}
L_n(\theta)\as L(\theta):=\frac{`E( |W|1_{\arg(W)\leq \theta})}{`E(|W|)},
\een
where $L(\theta)$ is the length of the curve $t\mapsto Z_\mu(t)$ between times $0$ and $F_\mu(\theta)$.
Fix a small $`e > 0$. There exists $\theta_1 < \cdots < \theta_k$ such that the convex hull of the points $Z_\nu(F_\nu(\theta_i))$ is at  distance at most $`e$ of $\BB_\nu$. Notice that such a property implies that the successive segments lengths $l_i=|Z_\nu(F_\nu(\theta_i))-Z_\nu(F_\nu(\theta_{i-1}))|$ satisfies
\[ L(\theta_{i})-L(\theta_{i-1})-2`e\leq l_i \leq L(\theta_{i})-L(\theta_{i-1})\]
since $B_{\nu}$ is convex and the graph of $Z_\nu$ must stay at distance at most $`e$ of
$[Z_\nu(F_\nu(\theta_i)),Z_\nu(F_\nu(\theta_{i-1}))]$ between times $F_\nu(\theta_{i})$ and $F_\nu(\theta_{i-1})$. But for $n$ large enough, up to an additional $\varepsilon$, the discrete curve has the same properties with high probability. By \eref{eq:deuz}
\[
\sup_{1\leq j\leq n} \left| \frac{S(N_n(\theta_j))}{n`E(|W|)} -
  Z_\nu(F_{\nu}(\theta_j)) \right|
\as 0.
\]
The length  $L_n(\theta_i)-L_n(\theta_{i-1})$ of the curve between $\theta_{i-1}$ and $\theta_i$ converges \AS to $L(\theta_{i})-L(\theta_{i-1})$ by \eref{eq:long}. This implies that the Hausdorff distance between $\bB_n/(n`E(|W|))$ and the convex hull of the points $\frac{S(N_n(\theta_j))}{n`E(|W|)}$'s goes to zero \AS  ~$\Box$

We now consider convolution and mixture of CCS.

\begin{prop}\label{pro:convol} Let $X$ and $Y$ be independent \rv 
  in $\Com$ with mean 0 (but not equal to 0 \AS), and
  $\lambda\in[0,1]$. Let $\mu_X$, $\mu_Y$ and $\mu_{X.Y}$ be the laws
  of $X$, $Y$ and $X.Y$. We have
\[\CC_{\sK(\mu_{X.Y})}=\CC_{\sK(\mu_X)}\star \CC_{\sK(\mu_Y)}\textrm{~~ and~~} \CC_{\sK(\lambda\mu_{X}+(1-\lambda)\mu_Y)}=\lambda\CC_{\sK(\mu_X)}+(1-\lambda) \CC_{\sK(\mu_Y)}.\]
\end{prop}

\Proof The statement concerning the mixture is quite easy and follows Theorem \ref{thm:cv-conve} for example. 
For the other one, following \eref{def:convo}, it suffices to see that
$\sK(\mu_{X.Y})=\sK(\mu_X)\CT \sK(\mu_Y).
$
Observe that for any measure $\mu$ on $\mathbb{C}$ (such that $0<|X_\mu|<+\infty$),
\[\frac{`E(e^{ix \arg(X_\mu)}|X_\mu|)}{`E(|X_\mu|)}=\int_{0}^{2\pi}e^{ix\theta}\frac{m_{X_\mu}(\theta)}{`E(|X_\mu|)} d`P_{\arg(X_\mu)}(\theta).\]
Indeed, according to \eref{eq:dnu}, the Fourier transform of $\sK(\mu)$ at position $x$ is given by $\frac{`E(e^{ix \arg(X_\mu)}|X_\mu|)}{`E(|X_\mu|)}$. 
Hence, the Fourier transform of $\sK(\mu_{X.Y})$, for $X$ and $Y$ independent, is 
\[\frac{`E(e^{ix \arg(X Y)}|X Y|)}{`E(|XY|)}=\frac{`E(e^{ix \arg(X)}|X|)}{`E(|X|)}\frac{`E(e^{ix \arg(Y)}|Y|)}{`E(|Y|)},\]
which implies that the Fourier transform of $\sK(\mu_{X.Y})$ and of $\sK(\mu_X)\CT \sK(\mu_Y)$ are the same. $\CC_{\sK(\mu_{X.Y})}$ and $\CC_{\sK(\mu_X)}\star \CC_{\sK(\mu_Y)}$ are equal by Definition \ref{def:convo}.~$\Box$

\begin{rem}
  \label{rem:noncharaccomplex} The CCS $\CC_{\sK(\mu)}$ characterises $\sK(\mu)$ but not $\mu$. 
For example the two following measures $\mu_1 =\frac{1}{3}\l( \delta(1) + \delta(e^{2i\pi/3}) +
    \delta(e^{4i\pi/3})\r)$ and  $\mu_2 =\frac{1}{3}\l( \frac{1}{2} \delta(\frac{1}{2}) + \frac{1}{2}
    \delta(\frac{3}{2}) + \delta(e^{2i\pi/3}) +
    \delta(e^{4i\pi/3})\r)$ satisfy  $\sK(\mu_1) = \sK(\mu_2)$ and $\CC_{\sK(\mu_i)}$ is  an equilateral triangle. Every CCS $\CC_\nu$ can therefore be  seen as an equivalence class of measures over $\Com$.\\
However, ${\sK\l(\mu_1 \CT \mu_1\r)}$ represents a polygon with $6$ sides, whereas ${\sK\l(\mu_1 \CT \mu_2\r)}$ a polygon with $7$ sides, even though $\sK(\mu_1) =
\sK(\mu_2)$. Hence $\sK(\mu_{1}\star \mu_2)$ is not a function of $\sK(\mu_1)$ and $\sK(\mu_2)$, and then the convolution of measures in $\mathbb{C}$ can not be turned into a nice operation on CCS.
\end{rem}

\section{Some models of random CCS}
\label{sec:mrc}
In this part, we consider the problem of finding natural distributions
on the set of CCS. We first recall some classical considerations on
simple models of random convex polygons.
In a second part we take advantage of the representation of CCS by
measures in $\MT^0$ to present models for the generation of smooth CCS
based on random Fourier coefficients.

\subsection{Reordering of closed polygons}
\label{sec:cccc}

Consider the problem of generating a convex polygon by specifying a
finite set of vectors representing its edges. 
Let $\mu$ be a distribution on $\mathbb{C}$ whose support is not
reduced to a point, and for some $n\geq 2$, let $(X_i,i=1,\dots,n)$ be
$n$ \iid \rv distributed according to $\mu$, and set
\[W_i=X_{(i \bmod{n}) +1}-X_i,~~1\leq i \leq n.\] %
Naturally, $\sum_{i=1}^n W_i=0$. Let
$(\hat W_i, 1\leq i\leq n)$ be the sequence $(W_i, 1\leq i \leq n)$
sorted according to their arguments. Let now $S$ be defined as in
\eref{eq:Sk}, and $\bB_n$ defined as in Section
\ref{sec:mes-C}. Further, let $\mu$ be the distribution of $W_1=X_2-X_1$,
and $\nu=\sK(\mu)$. 

\noindent The following result analogous with Theorem \ref{thm:cv-conve} shows 
that $\bB_n$ converges in distribution to $\BB_{\nu}$:

\begin{theo}\label{thm:conv-conv}%
   Assume that $\mu$ is centred (different from $\delta_0$). Then
   $$d_H\l(\bB_n/(n`E(|W|)), \BB_{\nu}\r)\as 0.$$
   Moreover \eref{eq:CVcompl} holds.
\end{theo}

\Proof 
We have $S(N_n(\theta))=\sum_{i=1}^n (X_{(i \bmod{n}) +1}-X_i )1_{\arg(X_{(i \bmod{n}) +1}-X_{i})\leq \theta}$; the difference with the proof of Theorem \ref{thm:cv-conve} is the dependence between the \rv in the sum. But these \rv are only weakly dependent (each \rv depends on the previous and following one); then strong law of large number applies to this case (since the sum can be split into two sums with \iid \rv), and the rest of the proof follows that of Theorem \ref{thm:cv-conve}.~$\Box$

\subsection{Convex polygon by conditioning / Convex polygon by chance}

Another natural way to sample a convex polygon is to take some
\iid points $W_0,\dots,W_{n-1}$ in the plane according to a distribution $\mu$ 
with support not included in a line, and to condition $(W_0,\dots,W_{n-1})$ to be a convex polygon. Define the set of all possible convex polygons as
\[\bB_n=\{{\bf w}:=(w_0,\dots,w_{n-1}) : \arg(w_{i+1 \bmod{n}}-w_i) \textrm{ forms an increasing sequence in }[0,2\pi)\}.\] %
Hence,  ${\bf w}$ represents the
list of vertices of a convex polygons encountered when following its boundary in the counter-clockwise direction (with some conditions for $w_0$).\par
The value of $\mu^{\otimes n}(\bB_n)$ is known only for $\mu$ equal to the uniform distribution in a triangle or in a parallelogram \cite{Valtr1,Valtr2} and in a circle \cite{MarckertCircle}; when $\mu$ is the uniform distribution in a CCS, the limit behaviour for $\bf w$ under the condition $\bf w\in \bB_n$ is described in \Barany \cite{BAR}. We open here a parenthesis to explain the underlying
difficulty. Consider $S_n:=(w_0,\dots,w_n)$ a $n$-tuple of points in
$\mathbb{`R}^2$, not three of them being on the same line (this
happens almost surely if $\mu$ admits a density on an open set in
$`R^2$). 
When $w_{i}=(x_i,y_i)$ for any $i$,  the algebraic area of the triangle $(w_i,w_j,w_k)$ is
\begin{equation}\label{eq:aire-tri}
A_{i,j,k}= \frac{1}{2}(x_iy_j+x_jy_k+x_ky_i-y_ix_j-y_jx_k-y_k x_i).
\end{equation}

The set $(s_{i,j,k}:={\sf sign}(A_{i,j,k}), 0\leq i<j<k\leq n-1)$ is called the
\textit{chirotope} of $S_n$. An equivalence class for the
chirotope, is called an \textit{order type}.  The sequence $S_n$ forms
a convex polygon iff all $s_{i,j,k}$ have the same sign. It is known that some
order types are empty, and also that deciding if an order type is not
empty, is a $NP$-complete problem (cf. Knuth \cite[Section 6]{Knuth}).
\par

When $(W_j=(X_j,Y_j),j=0,\dots,n-1)$ is a family of \iid \rv, such
that the $X_i$ and $Y_i$ are independent Gaussian centred \rv with
variance 1, it turns out that the Laplace transform of the joint law
of the $A_{i,j,k}$'s  (the areas of the triangles $(W_i,W_j,W_k)$) that is
\[\Phi(\lambda_{i,j,k}, 0\leq i<j<k\leq n-1):= `E\l(\exp\l(\sum_{0\leq i<j<k\leq n-1} \lambda_{i,j,k} A_{i,j,k}\r)\r)\]
is equal to $|\det(\Lambda)|^{-1/2}$, where $\Lambda = \left(
  \ell_{i,j}\right)$ and $\ell_{i,j} =
\sum_{a}{\lambda_{i,j,a}+\lambda_{a,i,j}-\lambda_{i,a,j}}$ (in a
neighbourhood of the origin of $`R^{\binom{n}{3}}$). To get this
result, the method is the same as the one for the computation of the
Fourier transform of a Gaussian vector in $\mathbb{R}^d$.

\begin{rem}
As remarked by Andrea Sportiello in a private communication,
$|\det(\Lambda)|$ is always a square of a polynomial in the
coefficients $\bar{\lambda}_{i,j}$. Indeed, for $\Lambda'=\begin{bmatrix}
           -Id_n & 0 \\
           0 &Id_n
\end{bmatrix}\Lambda$, $\Lambda$ and $\Lambda'$ have the same
determinant (up to factor $(-1)^n$).  But it can be shown
that $\Lambda'$ is a skew matrix, and then its determinant is
the square of its Pfaffian, which is indeed a polynomial on
its coefficients.
\end{rem} 

The Gaussian distribution is probably the simplest non trivial measure
for which this computation is possible. The question of the emptiness
of an order type $S=(s_{i,j,k},i<j<k)$ can be translated in term of
the support of the measure, but Knuth's result implies that it is a
difficult task. If $n=3$, only one triangle is present; the Laplace
transform is $1/(1-3\lambda_{0,1,2}^2/4)$, the transform of a Gamma
\rv with a random sign; when $n=4$, the Laplace transform is much more
complex.

\subsection{Generation of smooth random CCS}
\label{sec:positive-fourier-series}
This part is mainly prospective.
By Theorem \ref{th:mes-conv}, to conceive a model of random CCS in $\Conv(1)$ and to conceive a model of random measures with values in $\MT^0$ is the same problem. Since the condition ``to be in $\MT^0$'' has a simple expression in term of Fourier coefficients, and since the Fourier coefficients determine the measure (Proposition \ref{pro:mes-ser}), a simple idea consists in describing random measures in $\MT^0$ using random Fourier coefficients. \par
This leads us to  Szegö's Theorem \cite{SzegoOrthogonalPolynomials}: if a
trigonometric polynomial $P:\T\to `R^+$ admits only non-negative
values, then there exists a polynomial $D$ such that:
$$ \forall t \in \T, \qquad P(t) = |D(e^{it})|^2$$
Moreover $D$ is unique up to multiplication by a complex of
modulus $1$. If we consider the Fourier expansion
$D(e^{it})=\sum_{n\geq 0}\rho_ne^{i\theta_n} e^{i nt}$, for some
finite sequences of real numbers $(\rho_n), (\theta_n)$, the modulus
of $D$ is equal to:
\begin{align}
|D(e^{it})|^2&=A_0+\sum_{n\geq 1}A_n \cos(nt)+B_n\sin(nt) \notag
\\
\textrm{with}&\begin{cases}
A_0 =\sum_{k\geq 0} \rho_k^2\\
A_n =2\sum_{k\geq 0} \rho_{k+n}\rho_k \cos(\theta_k-\theta_{k+n}) 
\textrm{ for } n \geq 1, ~\\
B_n =2\sum_{k\geq 0} \rho_{k+n}\rho_k \sin(\theta_k-\theta_{k+n}) 
\textrm{ for } n \geq 1.
\end{cases}\label{eq:ABandrhotheta}
\end{align}

\noindent Hence, the trigonometric polynomial $P$ is the density of a
measure $\mu \in \MT^0$ iff the sequences $(A_n)$ and $(B_n)$ satisfy
$(i)$~the perimeter condition ($A_0=\frac{1}{2\pi}$,
ensuring that $\mu$ is a probability measure) and $(ii)$ the closed
path condition ($A_1=B_1=0$, ensuring that $\int_0^{2\pi}
e^{ix}d\mu(x)=0$).

\subsubsection{Generation of CCS via their Fourier coefficients}

In order to generate a random pair ${\cal P}:=((\rho_k,k\geq 0),
(\theta_k,k\geq 0))$ satisfying both conditions, two possibilities are
open, depending on which condition should be satisfied first (but the question of finding natural distributions for CCS will remain open). \par
To satisfy $A_1=B_1=0$ first, it
suffices to generate $\rho_{j}$ and~$\theta_j$ for $j\geq 1$ at
random then take $\rho_0$ and~$\theta_0$ such that:
\[\rho_0\rho_1 e^{i (\theta_0-\theta_1)}=-\sum_{k
  \geq 1} \rho_{k+1}\rho_k e^{i(\theta_k-\theta_{k+1})}.\]
This is always possible if the sum converges and if $\rho_1$ is not 0. To satisfy $A_0=1/2\pi$ from here, a normalisation step can be applied: divide each $\rho_n$ by $\sqrt{\sum_{k\geq 0}
  \rho_k^2}$. \par
\begin{figure}[htb]
  \begin{center}
  \begin{tabular}{ccc}
    \includegraphics[width=4.5cm]{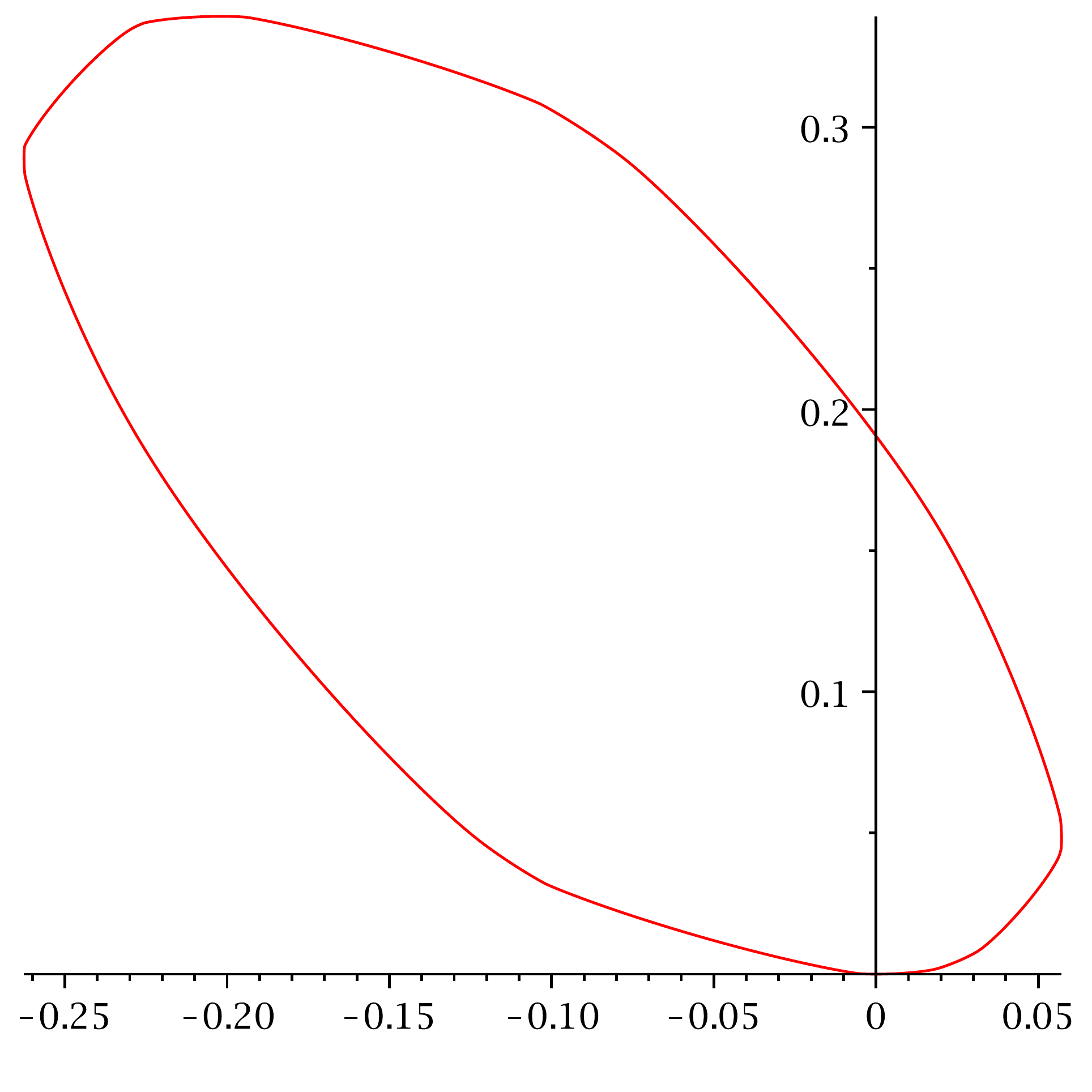} & 
    \includegraphics[width=4.5cm]{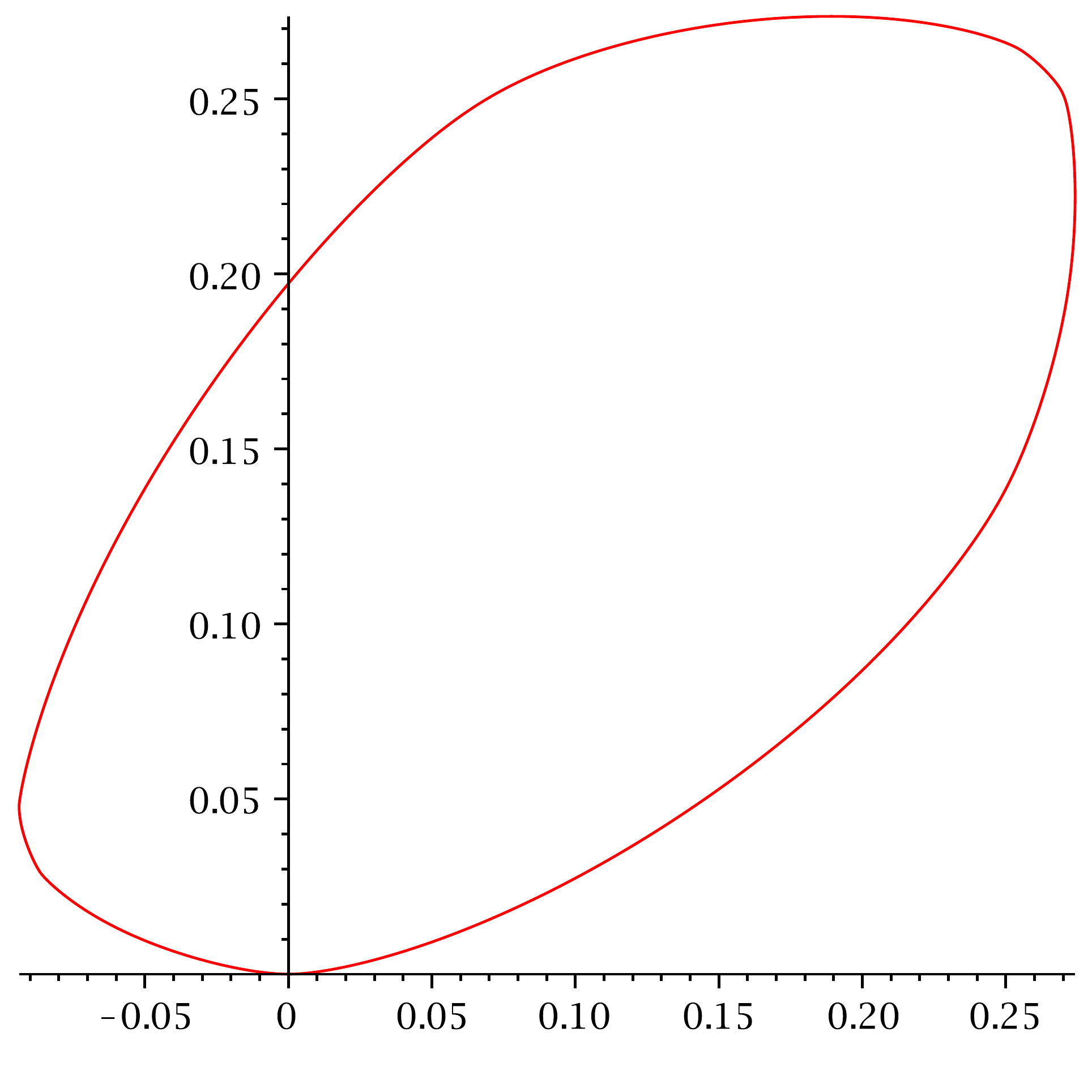} & 
    \includegraphics[width=4.5cm]{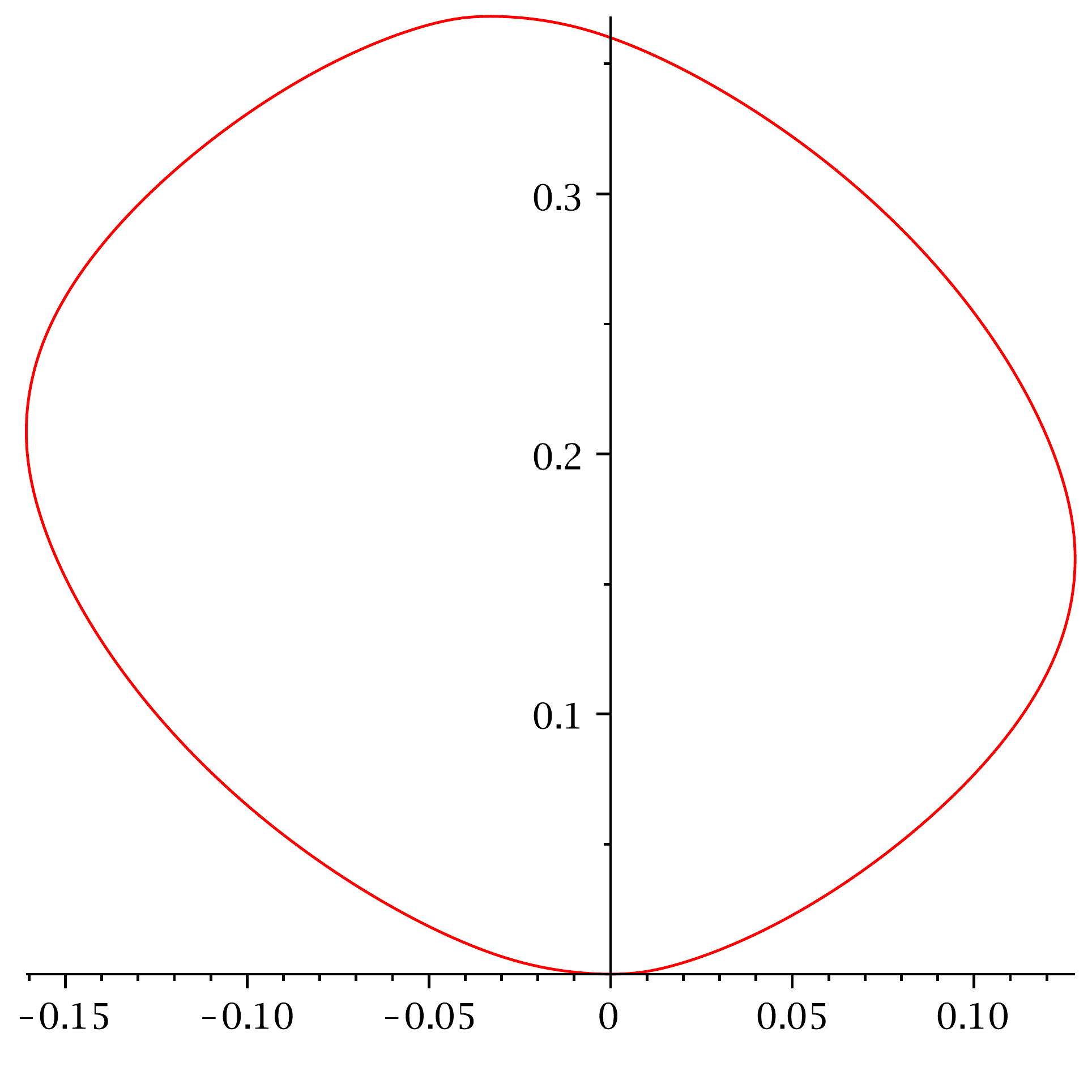} 
  \end{tabular}
  \end{center}
  \caption{Examples of random CCS generated from
    trigonometric polynomials containing 25 non-zero coefficients
    (with $\rho_j \sim \uniform [0;1]$, and
    $\theta_j \sim \uniform [0;2\pi]$), all these \rv being taken independently.}
  \label{fig:random_generation_1}
\end{figure}

Szegö's theorem ensures that the set of measures induced by this
method has full support over~$\MT^0$: indeed, each measures in
$\MT^0$ can be weakly approached by a sequence of distributions with
strictly positive density; these ones can be in turn approached by a
sequence of positive trigonometric polynomials, and Szegö's Theorem
gives a representation of these polynomials. The results of such a
generation can be seen on figure~\ref{fig:random_generation_1}.

Another solution consists in ensuring first $A_0=1/2\pi$, which comes
down to producing $(\rho_k,k\geq 0)$ such that $\sum_{k\geq 0}
\rho_k^2=\frac{1}{2\pi}$. This can be done by choosing (generating)
random reals $r_j$ in $[0,1]$, and setting:

\[\rho_k^2=\frac{1}{2\pi}r_k\prod_{j=0}^{k-1} (1-r_j).\] 
This is well defined if $\prod_k (1-r_k)$ converges to $0$ when $k$
goes to infinity (for example, taking \iid $r_j$'s under
$\uniform[0,1]$ does the job). From here, satisfying $A_1=0$ and
$B_1=0$ by a right choice of $\theta$'s can become more difficult, and
even impossible, for example if $\rho_0=\rho_1 > 0$ and all other
$\rho_i$'s are 0.
\begin{figure}[htb]
  \begin{center}
  \begin{tabular}{ccc}
    \includegraphics[width=4.5cm]{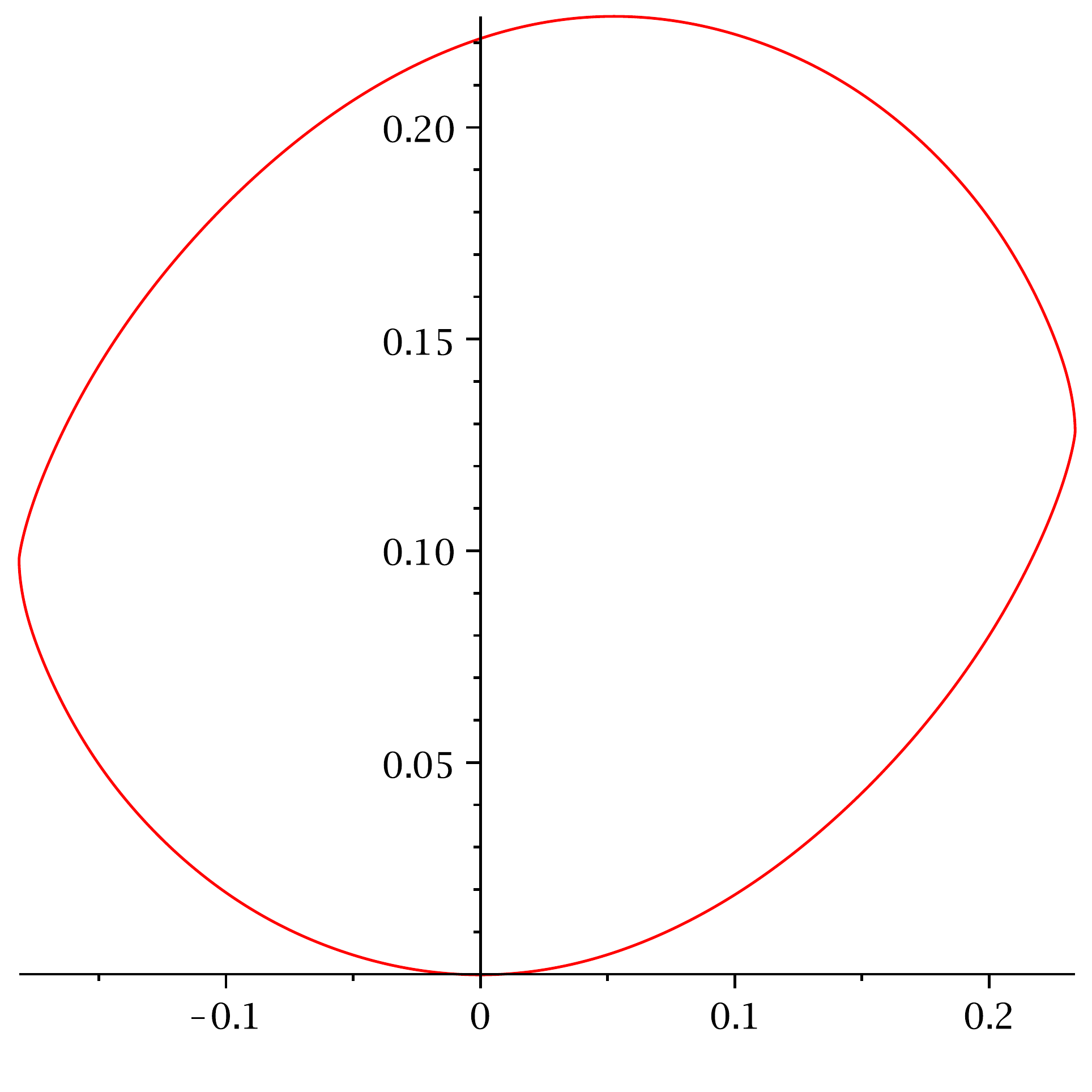} & 
    \includegraphics[width=4.5cm]{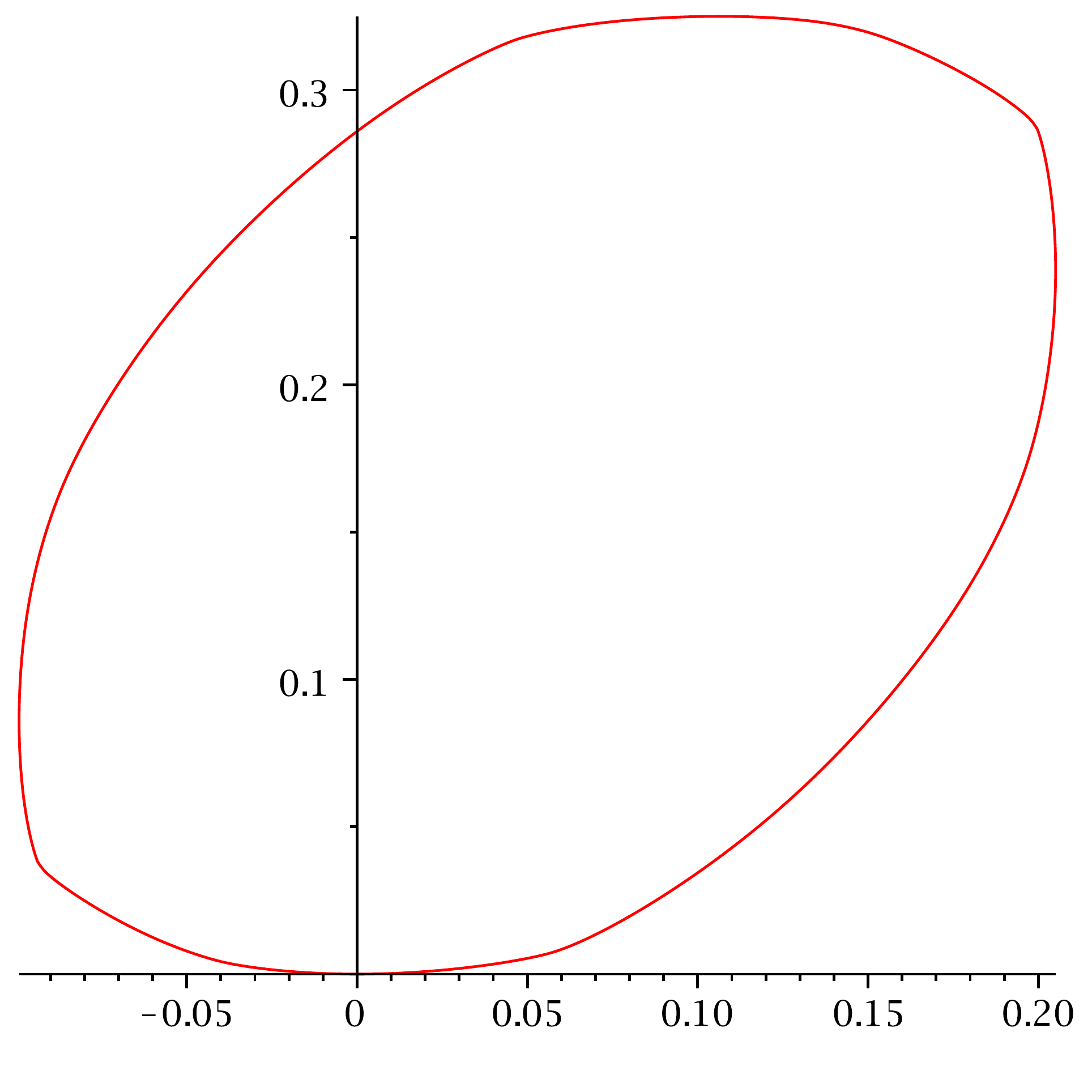} & 
    \includegraphics[width=4.5cm]{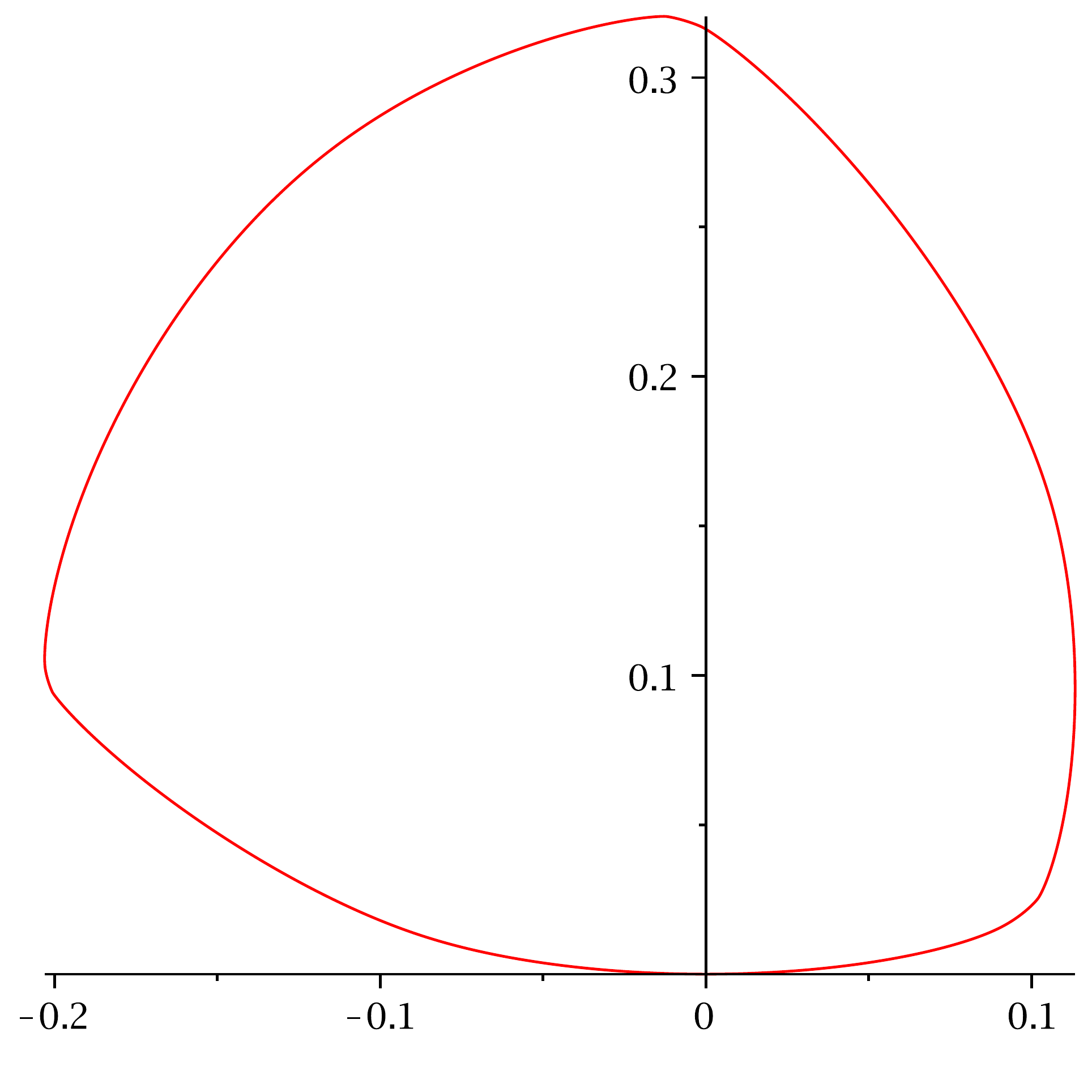} %\\
  \end{tabular}
  \end{center}
  \caption{Examples of random CCS generated from polynomials
    containing 12 non-zero coefficients with sparse coefficients (the
    indices of the non-null Fourier coefficients of $F$ are selected
    with probability $0$ if the previous coefficient was selected, and
    with probability $\frac{1}{2}$ otherwise; $\rho_j \sim \uniform
    [0;1]$; $\theta_j \sim \uniform [0;2\pi]$, all these \rv are
    taken independently).}
  \label{fig:random_generation_2}
\end{figure}
Nevertheless, it is possible to generate ${\cal P}$ satisfying all the
constraints at once. Choose (at random or not) a subset $F$ of
$\mathbb{N}$ such that if $i\in F$, then $i+1\notin F$, and a sequence
$x_k$ such that $\sum_{k\geq 0} x_k^2=\frac{1}{2\pi}$ as above.  Now,
let $n_j$ be the $j+1$-th smallest element in $F$, with the convention
that the smallest is $n_0$. Define the sequence $(\rho_k)$ by:
\[\rho_{n_j}=r_j, \quad \rho_k = 0~\textrm{otherwise}\]
\noindent Thanks to \eref{eq:ABandrhotheta}, $A_1 = B_1 = 0$ (since
for all $k$, $\rho_k\rho_{k+1} = 0$), and this for any choice of
$(\theta_k)$.  Examples of CCS generated this way appear on
Figure~\ref{fig:random_generation_2}.

\subsubsection{Generation of CCS with a given area}

\begin{figure}[tbh]
  \begin{center}
  \begin{tabular}{ccc}
    \includegraphics[width=4.5cm]{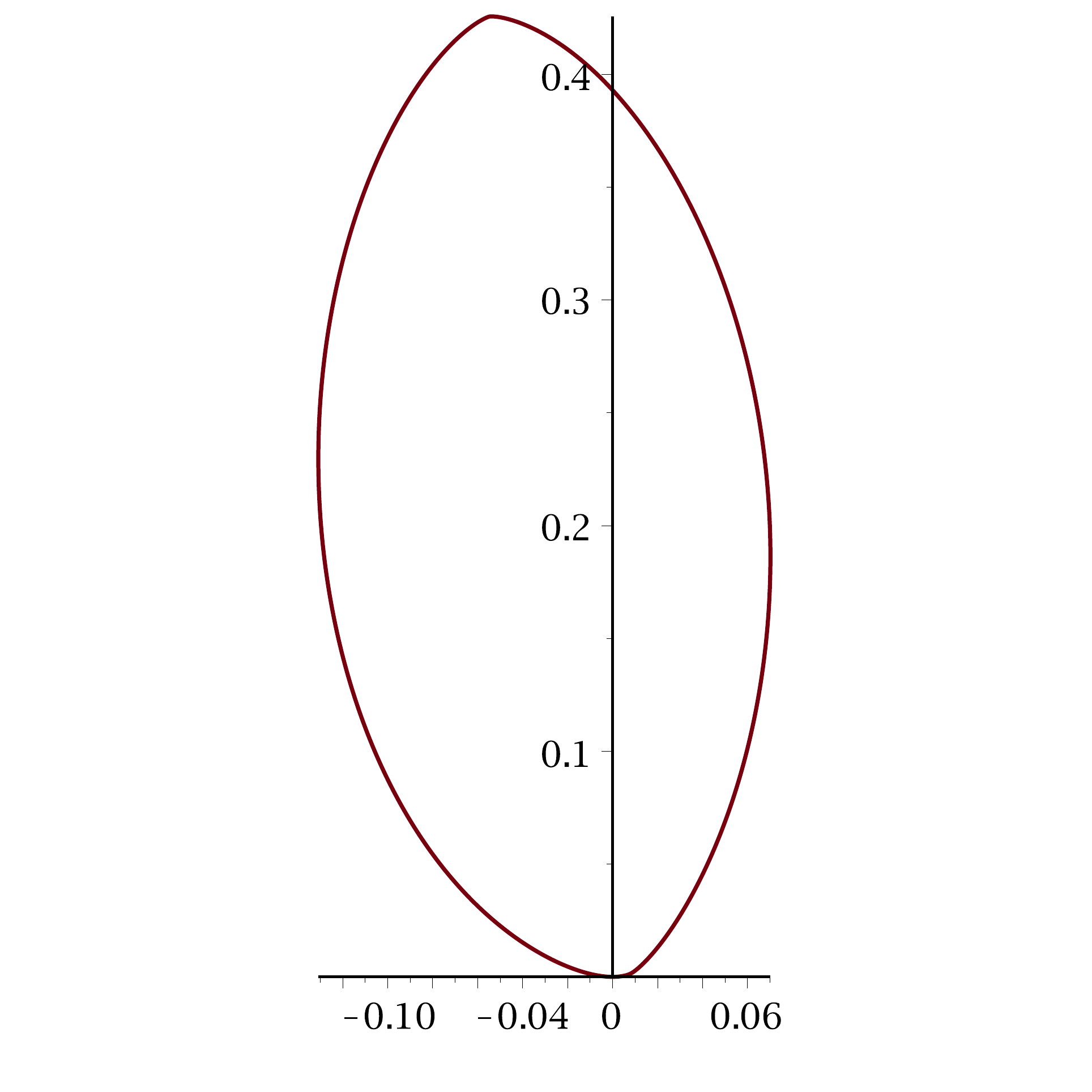} & 
    \includegraphics[width=4.5cm]{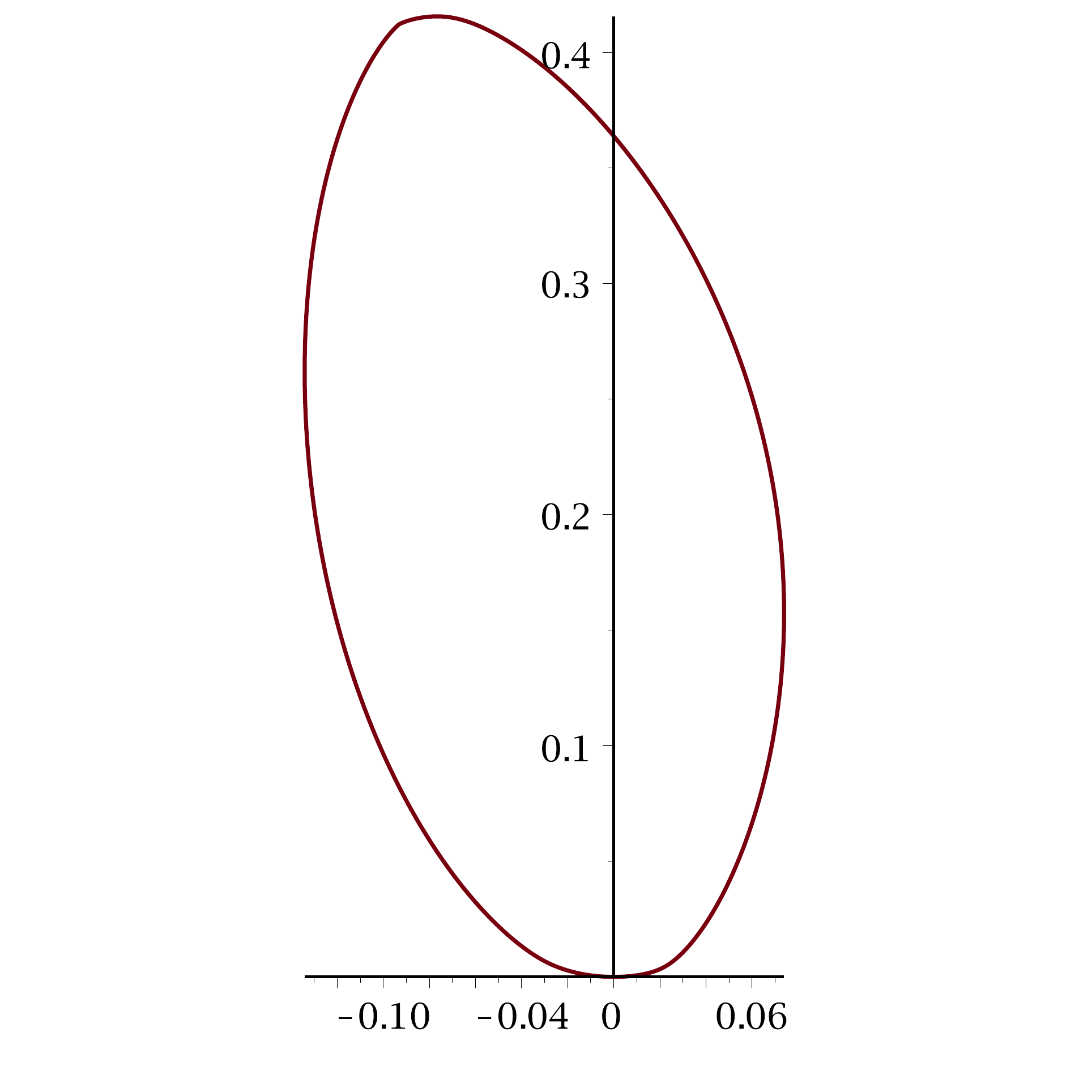} & 
    \includegraphics[width=4.5cm]{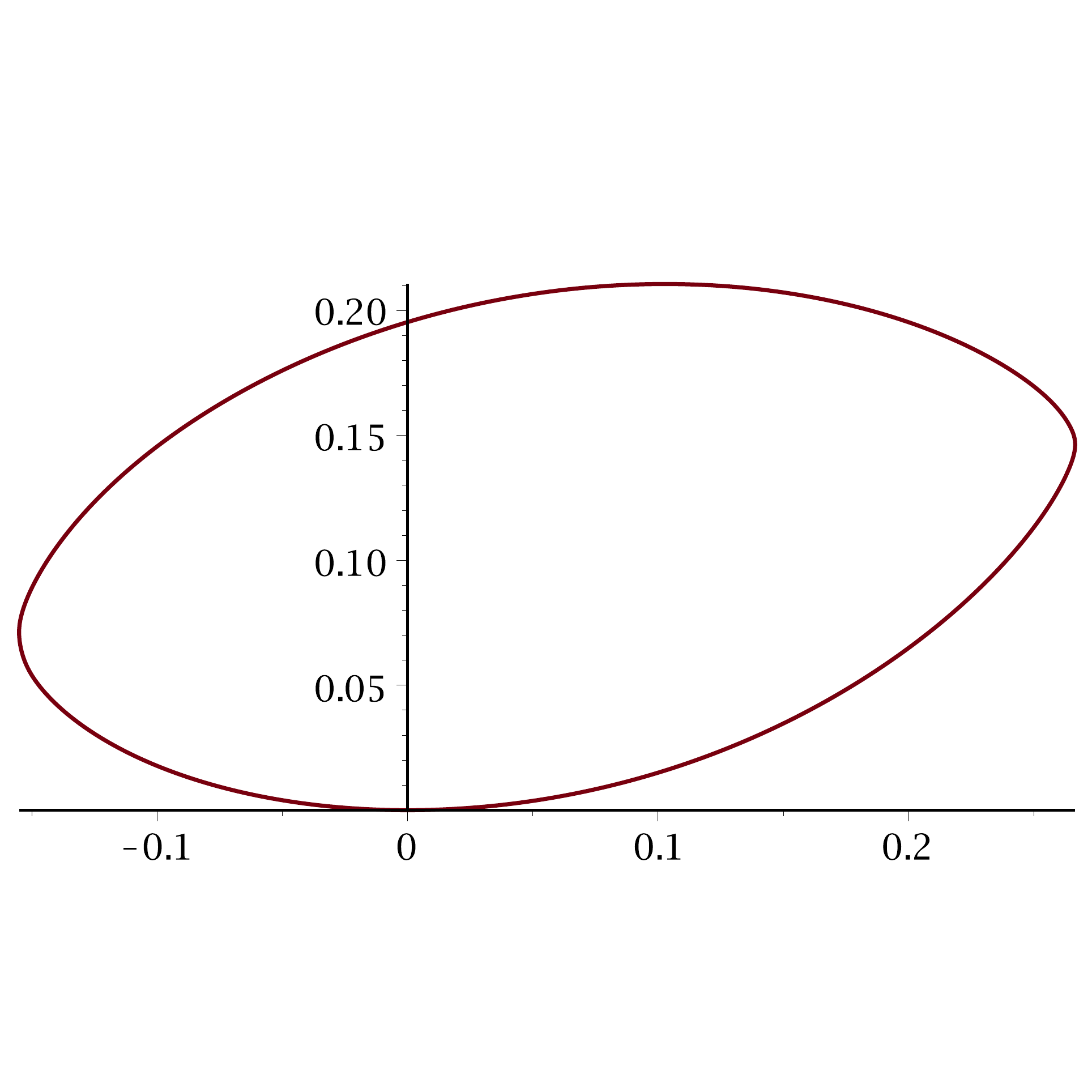} 
  \end{tabular}
  \end{center}
  \caption{Examples of random CCS of perimeter $1$ generated
    such that their area is equal to $\frac{1}{4\pi} - \frac{\pi}{2}
    \times 0.01$ (the polynomials possess 20 non-null coefficients,
    $\rho_j \sim \uniform [0;1]$, and $\theta_j \sim \uniform
    [0;2\pi]$, all these \rv being taken independently.)}
  \label{fig:random_generation_3}
\end{figure}

Consider the problem of generating a CCS in $\Conv(1)$ with a given area $\alpha=\frac{1}{4\pi} -\frac{\pi}{2} \beta \in [0, \frac{1}{2\pi^2}]$. Such a CCS corresponds to Fourier coefficients that satisfies:
\[\sum_{k\geq 2} \frac{a_k^2+b^2_k}{k^2-1}=\beta.\]
As in the previous section, we consider a sequence of
numbers $(r_j)$ in $[0,1)$ for $j\geq 2$, such that $\prod_{j\geq 2}
(1-r_j)=0$, and define positive reals $(c_k)$ such that:
\[ \frac{c_k^2}{k^2-1} = \beta ~r_k \prod_{j=2}^{k-1} (1-r_j).\]
Let $(\theta_k,k\geq 2)$ be a sequence of real numbers in $[0,2\pi)$. Then the Fourier coefficients of the associated measure can be computed as follows:
 \be a_k&=&\cos(\theta_k) c_k,~~
  b_k=\sin(\theta_k) c_k.  \ee
It is still possible to take $a_1=b_1=0$ and $a_0=1/(2\pi)$, but since we didn't use  Szegö's theorem, the standard Fourier series associated to the $a_i$'s and $b_i$'s is unlikely to be a
positive function. From here, it suffices to reject all series with a negative minimum.
The results of such a generation appear on
Figure~\ref{fig:random_generation_3}. Experiments show that the
rejection rate is very high, and that it is very difficult to generate
CCS with $\beta > 0.01$ (the theoretical maximum being 
$\frac{1}{2\pi^2}\approx 0.05$).

\section{Appendix}

\subsection{Proof of Theorem \ref{thm:cv-emp}}
\label{proof:thm:cv-emp}

\subsubsection*{Convergence of the FDD of $W_n$}
Let  $\theta_0:=0\leq \theta_1<\theta_2<\cdots<\theta_\kappa=2\pi$ for some $\kappa\geq 1$ be fixed. 
In the sequel, for any function (random or not) $L$ indexed by $\theta$, $\Delta L(\theta_j)$ will stand for $L(\theta_j)-L(\theta_{j-1})$.  For any $\ell\leq
\kappa$
\begin{equation}\label{eq:W}
W_n(\theta_\ell)=\sqrt{n}\sum_{j=0}^\ell \Delta \l[Z_n(N_n(\theta_j))-Z_\mu(F_\mu(\theta_j))\r]
\end{equation}
where  by convention $Z_n(N_n(\theta_{-1}))=Z_{\mu}(F_\mu(\theta_{-1}))=0$. The convergence of the FDD of $W_n$ follows from those of $(\sqrt{n}\Delta \l[Z_n(N_n(\theta_j))-Z_\mu(F_\mu(\theta_j))\r], 0\leq i \leq \kappa)$. Notice that 
\beq
\Delta Z_\mu(F_\mu(\theta_j))= `E(\exp(iX) 1_{\theta_{j-1}<X \leq \theta_j}).
\eq
If for some $j$, $\theta_{j-1}$ and $\theta_j$ are chosen in such a way that $\Delta F_\mu(\theta_j)=0$ then the $j$th increment in \eref{eq:W} is 0 almost surely (this is the case for the 0th increment if  $\mu(\{0\})=0$). We now discuss the asymptotic behaviour of the other increments : let $J=\{j \in \{0,\dots,\kappa\} : \Delta F_\mu(\theta_j)\neq 0\}$.

Let $(n_j,j \in J)$ be some fixed integers such that $n=\sum n_j$. Denote by $\mu_{\theta_{j-1},\theta_j}$ the law of $X_\mu$ conditioned on 
$\{\theta_{j-1}<X_\mu\leq \theta_j\}$, and by $X_{\theta_{j-1},\theta_j}$ a \rv under this distribution.  Conditionally on $(N_n(\theta_j)=n_j,j \in J)$, the \rv $\Delta Z_n(N_n(\theta_j)),j\in J$   are independent. The law of $\Delta Z_n(N_n(\theta_j))$ is that of a sum of $n_j-n_{j-1}$ \iid copies of \rv under $\mu_{\theta_{j-1},\theta_j}$, denoted from now on  $(X_{\theta_{j-1},\theta_j}(k),k\geq 1)$):
 \be
`E\l(\Delta Z_n(N_n(\theta_j)) \big| N_n(\theta_l)=n_l,l\in J\r)&=&n^{-1}`E\l(\sum_{m=1}^{n_j-n_{j-1}} e^{iX_{\theta_{j-1},\theta_j}(m)}\r)\\
&=&\frac{(n_j-n_{j-1})}n \frac{\Delta Z_\mu(F_\mu(\theta_j))}{\Delta
  F_\mu(\theta_j)}.  \ee 
Since $\l(\Delta N_n(\theta_j),j \in J\r)\sim \Mult\l(n, (\Delta F_\mu(\theta_{j}), j \in J)\r)$, 
\beq\label{eq:limit_N}
 \l(\frac{\Delta N_n(\theta_j)- n\Delta F_\mu(\theta_{j})}{\sqrt{n}} ,j \in J\r)\dd (N_j,j\in J)\eq
where $(N_j,j \in J)$ is a centred Gaussian vector with covariance function
\[\cov(N_k,N_l)=-\Delta F_\mu(\theta_{k})\,.\, \Delta F_\mu(\theta_{l}),\]
formula valid for any $0\leq k,l\leq \kappa$. Putting together the previous considerations, we have, conditioning first on the $N_n(\theta_j)$'s, and then integrating on the distribution of these \rv,
\ben\label{eq:deltaW_n}
\Delta W_n(\theta_j)= \sum_{l=1}^{\Delta N_n(\theta_j)} \frac{e^{iX_{\theta_{j-1},\theta_j}(l)}-`E(e^{iX_{\theta_{j-1},\theta_j}})}{\sqrt{n}}+\l(  \frac{\Delta N_n(\theta_j)- n\Delta F_\mu(\theta_j)}{\sqrt{n}}\r) `E(e^{iX_{\theta_{j-1},\theta_j}})
\een
Using \eref{eq:limit_N} and the central limit theorem, we then get that
\beq\label{eq:FDD}
(\pi \Delta W_n(\theta_j),0\leq j \leq \kappa)\dd  \sqrt{\Delta F_\mu(\theta_j)}\widetilde{N}_j+N_j
\begin{bmatrix}
`E(\cos(X_{\theta_{j-1},\theta_j}))\\
`E(\sin(X_{\theta_{j-1},\theta_j}) 
\end{bmatrix},
 \eq
where the \rv $N_j,\tilde{N_j},j\leq \kappa$ are independent, and the \rv $\tilde{N_j}$ are centred Gaussian \rv with covariance matrix, the covariance matrix of $\begin{bmatrix}
\cos(X_{\theta_{j-1},\theta_j})\\
\sin(X_{\theta_{j-1},\theta_j}) 
\end{bmatrix}$.

\subsubsection*{Tightness of $\{W_n, n \geq 0\}$ in $D[0,2\pi]$}

A criterion for tightness in $D[0,2\pi]$ can be found in Billingsley \cite[Thm. 13.2]{BIL}: a sequence of processes $(W_n,n\geq 1)$ with values in $D[0,2\pi]$ is tight if, for any $`e\in(0,1)$, there exists $\delta>0, N >0$ such that
 \[\lim_{\delta\to 0} \limsup_n`P(`o'(W_n,\delta)\geq `e)=0\]
where $`o'(f,\delta)=\inf_{(t_i)}\max_i \sup_{s,t\in[t_{i-1},t_i)} |f(s)-f(t)|,$ and the partitions $(t_i)$ range over all partitions of the form $0=t_0<t_1<\cdots<t_n\leq 2\pi$ with $\min\{t_i-t_{i-1}, 1\leq i\leq n\}\geq \delta$. 

 Since only the tightness in $D[0,2\pi]$ interests us, we will focus
 on $\Re(W)$ (since the imaginary part can be treated likewise, and
 since the tightnesses of both $\Re(W)$ and $\Im(W)$ implies that
 of~$W$). For the sake of brevity, in the sequel, we will use $W$ instead of $\Re(W)$.\par

The first step in our proof consists in comparing the distribution $`P_n$ of a set $\{X_1,\dots,X_n\}$ of $n$ \iid copies of $X_\mu$ with a Poisson point process $P_n$ on $[0,2\pi]$ with intensity $n \mu$, denoted by $`P_{P_n}$.
Conditionally on $\#P_n=k$, the $k$ points $P_n:=\{Y_1,\dots,Y_k\}$ are \iid and have distribution $\mu$, and then
$`P_{P_n}(~\cdot~ |\# P=n)=`P_n$. The Poisson point process is naturally equipped with a filtration $\sigma:=\l\{\sigma_t=\sigma(\{ P\cap[0,t]\}),t\in[0,2\pi]\r\}$.

We are here working under $`P_{P_n}$, and we let $N(\theta)=\# P_n \cap [0,\theta]$; notice that under $`P_n$, $N$ and $N_n$ coincide.\par

We will show the tightness of $W$ under $`P_{P_n}$ first. Before doing
this, let us see why it implies the same result under $`P_n$: Let
$m$ be a point in $[0,2\pi]$ such that $F_\mu(x)> 1/4, 1-F_\mu(x)>1/4$ (it is a kind of median of $\mu$). We need in the sequel $1-F_\mu(m)>0$; for measures in $\MT^0$ this is always the case, since if not, an atom with weight $>1/2$ would exist. We
will see that the tightness under $`P_{P_n}$ implies that the sequence
of processes $W$  under $`P_n$ is tight in $D[0,m]$ (the same proof works on $D[m,2\pi]$ by a time reversing argument). We claim that for any event $\sigma_m$ measurable,
\beq\label{eq:inter}
`P_n(A)=`P_{P_n}(A \,|\, \#P=n) \leq c\, `P_{P_n}(A)
\eq 
for a constant $c$ independent on $n$ and of $A$ (but which depends on $\mu$). 
This in hand, the tightness under $`P_{P_n}$ of $W$ on $D[0,m]$ implies that under $`P_n$.
Let us prove \eref{eq:inter}. We have
\be
`P_{P_n}(A \,|\, \#P=n)
         &=&\sum_{k} \frac{`P_{P_n}(A ,\#(P\cap[0,m])=k)`P(\#P\cap[m,2\pi]=n-k)}{`P(\#P=n)}\\
         & \leq &\sum_{k} `P_{P_n}(A ,\#(P\cap[0,m])=k)\sup_{k'}\frac{`P(\#P\cap[m,2\pi]=n-k')}{`P(\#P=n)}\\
         & \leq & c\,`P_{P_n}(A)
\ee
where $c=\sup_{n\geq 1}
\sup_{k'}\frac{`P(\#P\cap[m,2\pi]=n-k')}{`P(\#P=n)}$, which is indeed
finite since:

\begin{itemize}
\item first $\#P\cap[m,2\pi]\sim \Pois(n (1-F_\mu(m)))$, and then
  $\sup_{k'}`P(\#P\cap[m,2\pi]=n-k')$ is the mode of a Poisson
  distribution. When the parameter is $\lambda$, the mode is
  equivalent to $1/\sqrt{2\pi \lambda}$ when $\lambda\to+\infty$, so
  here it is equivalent to $1/\sqrt{2\pi n (1-F_\mu(m)) }$,
\item and by Stirling $`P(\#P=n)\sim (2\pi n )^{-1/2}$.
\end{itemize}

Working with a Poisson point process instead of working with $n$ \rv provides some independence between the number of \rv $X_i$ in disjoint intervals, and then on the fluctuations of $W_n$ in disjoint intervals.

Before starting, recall that if $N\sim \Pois(a)$, for any positive $\lambda$,
 \ben\label{eq:c-p}
`P(N\geq x) &=&`P(e^{\lambda N}\geq e^{\lambda x})\leq `E(e^{\lambda N-\lambda x})
=e^{-a+ae^{\lambda}-\lambda x}\\
`P(N\leq x)&=&`P(e^{-\lambda N}\geq e^{-\lambda x}) \leq `E(e^{-\lambda N+\lambda x})=e^{-a+ae^{-\lambda}+\lambda x}.
\een 
Let $A_\mu=\{x \in [0,2\pi], \mu(\{x\})>0\}$ be the set of positions of the atoms of $\mu$.
We now decompose $\mu=\mu|_{A_\mu} +\mu|_{\complement A_\mu}$; under $`P_n$ as well as under $`P_{P_n}$, the process $W$ can be also decomposed under the form $W|_{A_\mu}+W|_{\complement A_\mu}$  using $N|_{A_\mu}(\theta)=\# P\cap [0,\theta]\cap A_\mu$, $Z|_{A_\mu}(N|_{A_\mu}(\theta))=\sum_{j=1}^{N} e^{i\hat{X}_j}1_{\hat{X}_j\in A_\mu}$, etc. The 
fluctuations of $W=W|_{A_\mu}+W|_{\complement A_\mu}$ are then bounded by the sum of the fluctuations of both processes $W|_{A_\mu} $ and $W|_{\complement A_\mu}$. It is then sufficient to show the tightness for a purely atomic measure $\mu$, and for a measure having no atom $\mu$.

\subsubsection*{Case where $\mu$ is purely atomic}

Take some (small) $\eta\in(0,1)$, $`e>0$; we will show that one can find a finite partition $(t_{i},i \in I)$ of $[0,2\pi]$ and a  $\delta\in(0,1)$ such that 
\beq\label{eq:tightness}
\limsup_n`P_n(`o'(W_n,\delta)\geq `e)\leq \eta,
\eq
which is sufficient for our purpose. In fact we will establish \eref{eq:tightness} under $`P_{P_n}$ instead, on $[0,m]$ and then on $[m,2\pi]$, since we saw that this was sufficient (replacing $\eta$ by $c\eta$ in \eref{eq:tightness}, suffices too). 

Now, let $A_\mu^{\ge a}:=\{x \in A_\mu : \mu(\{x\})\geq a\}$. Clearly
$\#A_\mu^{\ge a}\leq 1/a$ and $[0,2\pi]\setminus A_\mu^{\geq a}$ forms
a finite union of open connected intervals $(O_{x},x\in G)$, with
extremities $(t'_i,i \in I)$. The intervals $(O_x,x\in G)$ can be
further cut as follows:

\begin{itemize}
\item do nothing to those such that $\mu(O_x)<2a$,
\item those such that $\mu(O_x)>2a$ are further split. Since they contain no atom with mass $>a$, they can be split into smaller intervals having all their weights in $[a,2a]$ except for at most one (in each interval $O_x$ which may have a weight smaller than $a$).
\end{itemize}

Once all these splittings have been done, a list of at most $3/a$ intervals are obtained, all of them having a weight smaller than $2a$. Name $G_a=(O_x, x\in I_a)$ the collection of obtained open intervals, index by $I_a$, and by $(t_i^a, i\geq 0)$ the partitions obtained. Clearly
\[M_a:=\max_{i\in I_a} `E(\cos(X_\mu)^21_{X_\mu \in O_i})\leq M'_{a}:=2a.\]
\subsubsection*{Control of the fluctuations of $W_n$ on an interval $O_x$}

In the sequel we take $a=`e^3$ and consider a unique interval $O_x=(\theta_{j-1},\theta_j)\in G_a$, in which case
we have $M_{`e^3}\leq 2`e^3$. We control first the last position of the random walk $W_n$. Under $`P_{P_n}$, ${\cal P}(n\mu\{\theta\}):=\#P_n\cap \{\theta\}$ has distribution $\Pois(n\mu(\{\theta\}))$, the \rv corresponding to different points being independent. Following \eref{eq:deltaW_n}, under $`P_{P_n}$, we get
\beq\label{eq:contro}
\Delta W_n(\theta_j)= \sqrt{n}   \sum_{
  \jfmatop{\theta \in A_\mu}{\theta_{j-1}\leq  \theta<\theta_j} } 
\l( \frac{{\cal P}(n\mu\{\theta\})}n- \mu(\{\theta\})\r) \cos(\theta).
\eq
These centred \rv can be controlled as usual Poisson \rv as recalled above. On the first hand, 
\ben\label{eq:prem2}
`P(\Delta W_n(\theta_j)\geq `e)&=& `P\l( \sum_{\theta }{\cal P}(n\mu\{\theta\})\cos(\theta)\geq y\r)
\een
where 
\beq
y=`e\sqrt{n}+n`E(\cos(X)1_{X\in A_\mu, \theta_{j-1}< X\leq \theta_j})
\eq
and where the set of summation is the same as before (from now on, it will
be omitted). Writing 
$`P\l( \sum_{\theta }{\cal P}(n\mu\{\theta\})\cos(\theta)\geq y\r)\leq \inf_{\lambda >0} e^{-\lambda y}\prod_{\theta}`E(e^{(\lambda\cos(\theta)){\cal P}(n\mu\{\theta\})})$
one has
\be `P(\Delta W_n(\theta_j)\geq `e)&\leq & \inf_{\lambda>0} \exp\l(-\sum_{\theta }n\mu\{\theta\}+\sum_{\theta }n\mu\{\theta\}e^{\lambda\cos(\theta)}   -\lambda y\r) .
\ee
To get a bound we will take $\lambda= {`e}/(2{\sqrt{n}M'_{`e^3}})$. This allows one to bound $e^{\lambda\cos(\theta)}$ by $1+\lambda\cos(\theta)+\lambda^2 \cos(\theta)^2$ which is valid uniformly for any $\theta$ provided that $n$ is large enough. Hence for $n$ large enough replacing $y$ by its value,
\be 
`P(\Delta W_n(\theta_j)\geq `e)&\leq & \inf_{\lambda>0} 
\exp\l(\lambda^2 n`E(\cos^2(\theta)1_{\theta\in I_x})   -\lambda `e\sqrt{n}\r) \\
&\leq & \inf_{\lambda>0} 
\exp\l(\lambda^2 n M'_{`e^3}  -\lambda `e\sqrt{n}\r) \\
&\leq & \exp(-1/(4`e))
\ee
this last equality being obtained for $\lambda=`e/(2M'_{`e^3}\sqrt{n})$.

The proof for the control of $`P(\Delta W_n(\theta_j)\leq -`e)\leq \inf_{\lambda>0} `E\l(e^{-\lambda \Delta W_n(\theta_j)-\lambda\delta}\r)$  for $\delta>0$ gives rise to the same estimates, except that the bound  $e^{\lambda\cos(\theta)}$ by $1-\lambda\cos(\theta)+\lambda^2 \cos(\theta)^2/4$ is taken to replace the other one, giving a bound  $\exp(-1/(2`e))$ at the end. 

Now we have to control the fluctuations, and not only the terminal value of the random walk. Theorem 12 p.50 in Petrov \cite{PET} allows one to control the first ones using the second ones.

\subsubsection*{Control of the fluctuations of $W_n$ on all  intervals}
The control of all intervals all together can be achieved using the union bound : since they are at most $3/`e^3$ such intervals by the union bound 
\[`P_{P_n}(\sup_{j} \Delta W_n(\theta_j)\geq `e)\leq {3}{`e^{-3}}e^{-1/(4`e)}.\]
This indeed goes to 0 when $`e\to 0$. 

\subsubsection*{Case where $\mu$ has no atom}
We now show the tightness of $W$ under $`P_{P_n}$ when $\mu$ has no atom and use the same method as before: we work under $`P_{P_n}$, cut $[0,2\pi]$ under sub-intervals $[t_{j-1},t_j]'$s, control the differences between starting and ending values on these intervals, since we saw that it was sufficient. \par 
First we cut  $[0,2\pi]$ into  $n$ (tiny) equal parts $([2\pi (j-1)/n,2\pi j/n],j=1,\dots,n)$. From \eref {eq:deltaW_n}
\beq\label{eq:diff-W}
 W(2\pi j/n)-W(2\pi j'/n)= \sum_{l=j'+1}^j \Gamma_l +\Theta_l
\eq
where, under $`P_{P_n}$, denoting further $\theta_j=2\pi j/n$,
\be
\Gamma_l & = & \sum_{m=1}^{{\cal P}(n \Delta(F_\mu(\theta_{l})))} \frac{\cos(X_{\theta_{j-1},\theta_j}(m))-`E(\cos(X_{\theta_{j-1},\theta_j}))}{\sqrt{n}}\\
\Theta_l & = & \frac{ {\cal P}(n \Delta(F_\mu(\theta_{l})))- n\Delta F_\mu(\theta_l)}{\sqrt{n}} `E(\cos(X_{\theta_{l-1},\theta_l}))
\ee
and ${\cal P}(\lambda)\sim \Pois(\lambda)$ and the different Poisson \rv appearing in the $\Gamma_l$ and $\Theta_l$ are independent. Let $`e>0$ be given and $N_{`e^3}=\lceil 1/`e^3\rceil$.  
Since $\mu$ has no atom there exists some times $t_0=0<t_1,\dots<t_{N_{`e}}=2\pi$ such that $\mu([t_{i-1},t_i])\leq `e^3$. 
We now control the fluctuations of $W$ on these intervals.

Write $D_j:=W(\frac{\lfloor 2\pi t_j n \rfloor}{n})-W(\frac{\lfloor 2\pi t_{j-1} n\rfloor}{n})$ as a sum of \rv $\Gamma_l$ and $\Theta_l$ as in \eref{eq:diff-W}:
\[D_j=S_j+S_j'\]
where\[S_j=\sum_{l=\lfloor 2\pi t_{j-1} n\rfloor +1}^{\lfloor2\pi t_{j}n\rfloor} \Gamma_l, \qquad S'_j= \sum_{l=\lfloor 2\pi t_{j-1} n\rfloor+1}^{\lfloor2\pi t_{j}n\rfloor}\Theta_l.\]
Each $\Gamma_l$ is itself a sum which involves a Poisson number of terms: the total number of terms in $S_j$ is $N_{t_j}-N_{t_{j-1}}$, a Poisson \rv  with parameter smaller than $`e^3 n$ under $`P_{P_n}$. From \eref{eq:c-p},   $`P_{P_n}(N(t_j)-N(t_{j-1})\geq 3`e^3n)\leq e^{-c`e^3n}$ for some positive $c$, this meaning that with high probability, $S_j$ is a sum of less than $3`e^3n$ centred and bounded \rv of the form  $\frac{\cos(X_{\theta_{j-1},\theta_j}(m))-`E(\cos(X_{\theta_{j-1},\theta_j}))}{\sqrt{n}}$. By Hoeffding's inequality 
\[`P(|S_j|\geq `e | N(t_j)-N(t_{j-1})\leq 3`e^3n)\leq c' \exp(-c/`e)\] for some $c,c'>0$.

The sum $S'_j$ is controlled as above, in the atomic case (see \eref{eq:contro} and below).~\\

We now show 2); since $f\mapsto \max_\theta |f(\theta)|$ is continuous on $D[0,2\pi]$, we only need to prove  $d_H(B_n,\BB_\mu)=\max_\theta |Z_n(N_n(\theta)/n)-Z_\mu(F_\mu(\theta))|$.

Since $B_n$ and $\BB_\mu$ are compact, there exists $(x_n,x)$ in $B_n\times \BB_\mu$ realising this distance: $|x_n-x|=d(x_n,\BB_\mu)=d(B_n,x)=d_H(B_n,\BB_\mu)$. Consider now the set of directions $\Theta_n$ and $\Theta$ of the tangents at $x_n$ on $B_n$ and that at $x$ on $\BB_\mu$ (we call here a tangent at $a$ on $A$ a line $l$ that passes by $a$ and such that $A$ is contained in one of the close half plane defined by $l$. The set of directions of these tangents is an interval). We claim that there exists in $\Theta_n\cap \Theta$ the direction $\theta^\star$ orthogonal to $(x_n,x)$. If not, this means that at $x_n$ (or at $x$) the line passing at $x_n$ (or $x$) and orthogonal to $(x_n,x)$ crosses $B_n$ (or  $\BB_\mu$). This would imply that in a neighbourhood of $x$ (or $x_n$) there exists a point $x'$ (or $x'_n$) closer to $x_n$ (resp. $x$) than $x$ (resp. $x_n$), a contradiction. 

To end the proof, we need to show that $(x, x')$ corresponds to some $(S_n(N_n(\theta)/n),Z_\mu(F_\mu(\theta)))$. In other words, they are extremal points on their respective curves, and owns some parallel tangents. The second statement is clear.
For the first one, we have to deal with the fact that $B_n$ (and so do $\BB_{\mu}$ for certain measures $\mu$) have linear portions. But the distance between $B_n$ and $\BB_{\mu}$ is not reached inside the linear intervals since the Hausdorff distance between a segment $[a,b]$ and a CCS $C$ is given by $\max\{d(a,C),d(b,C)\}$.   
~$\Box$

\subsection*{Acknowledgements} 
We thank both referees for their numerous remarks that really helped to improve the paper.
\small
\bibliographystyle{plain}
%\bibliography{convex}

%\newpage
%\small
\setcounter{tocdepth}{3}
%\tableofcontents

\end{document}